\documentclass[11pt,reqno,oneside]{amsart}
\usepackage{amscd}
\usepackage{amsmath}
\usepackage{latexsym}
\usepackage{amsfonts}
\usepackage{amssymb}
\usepackage{amsthm}
\usepackage{graphicx}
\usepackage{verbatim}
\usepackage{mathrsfs}
\usepackage{enumerate}
\usepackage{hyperref}
\usepackage{fullpage}
\usepackage{xcolor}
\usepackage[latin1]{inputenc}

\theoremstyle{plain}
\theoremstyle{plain}\newtheorem{theorem}{Theorem}[section]
\theoremstyle{plain}\newtheorem{lemma}[theorem]{Lemma}
\theoremstyle{plain}\newtheorem{coro}[theorem]{Corollary}
\theoremstyle{plain}\newtheorem{proposition}[theorem]{Proposition}
\theoremstyle{plain}\newtheorem{remark}{Remark}[section]

\numberwithin{equation}{section}

\newcommand{\norm}[1]{\left\|#1\right\|}

\newcommand{\R}{\mathbb{R}}
\newcommand{\be}{\begin{equation}}
\newcommand{\ee}{\end{equation}}
 \newcommand{\ba}{\begin{aligned}}
 \newcommand{\ea}{\end{aligned}}

  \newcommand{\ben}{\begin{enumerate}}
   \newcommand{\een}{\end{enumerate}}

\newcommand{\Rmnum}[1]{\expandafter\@slowromancap\romannumeral #1@}

\allowdisplaybreaks

\makeatletter
\@namedef{subjclassname@2020}{\textup{2020} Mathematics Subject Classification}
\makeatletter

\begin{document}
\title{On the blowup of solutions for a nonlocal multi-dimensional transport equation}
\author{Wanwan Zhang$^{1*}$}

\address{$^1$ School of Mathematics and Statistics, Jiangxi Normal University, Nanchang, 330022, Jiangxi,
P. R. China}
\email{zhangww@jxnu.edu.cn}

\subjclass[2020]{35Q35; 35B44; 76D05}
\keywords{Transport equation; Nonlocal velocity; Blow-up}
\thanks{$^*$Corresponding author}

\begin{abstract}
In this paper, we revisit the problem of finite-time blowup for a multi-dimensional nonlocal transport equation studied in [Dong, Adv. Math. 264
(2014) 747-761]. Inspired by a one-dimensional analogous model considered in [Li-Rodrigo, Adv. Math. 374 (2020) 107345], we establish a new weighted nonlinear inequality implying the blow-up by a completely real variable based technique.
In particular, an inequality for the Riesz transform is obtained.
\end{abstract}
\smallskip
\maketitle
\section{Introduction and main results}
This paper is concerned with a multi-dimensional transport equation with a nonlocal velocity
\begin{equation}\label{M-CCF-T}
\left\{\ba
&\partial_t\theta+u\cdot \nabla\theta= 0, ~(x,t)\in \R^n\times\R_+,\\
&u=\nabla\Lambda^{-2+2\alpha}\theta,\\
&\theta(x,0)=\theta_{0}(x),~x\in\R^n, \ea\ \right.
\end{equation}
 where $n\geq2$ and $0<\alpha<1$.  Here the unknown $\theta$ defined in $\R^n\times\R_+$ is a scalar function, and the fractional Laplacian $\Lambda^s=(-\Delta)^{\frac{s}{2}}$ with $s\in \R$ is defined through the Fourier transform \cite{[Stein]}:
 \begin{eqnarray*}
\widehat{\Lambda^sf}(\xi)=(2\pi|\xi|)^s\widehat{f}(\xi).
\end{eqnarray*}
The second equation in \eqref{M-CCF-T} makes this model nonlocal.
The fractionally dissipative version of \eqref{M-CCF-T} reads as
\begin{equation}\label{dissipative}
\left\{\ba
&\partial_t\theta+u\cdot\nabla \theta+\Lambda^\gamma\theta = 0, ~(x,t)\in \R^{n}\times\R^{+}\\&u=\nabla\Lambda^{-2+2\alpha}\theta\\
&\theta(x,0)=\theta_{0}(x),~x\in\R^{n} \ea\ \right.
\end{equation}
where $n\geq1$, $0<\alpha<1$ and $0<\gamma<2$.
\par
The nonlocal active scalar equations \eqref{M-CCF-T} and \eqref{dissipative} were extensively studied.
When $n=1$ and $\alpha=\frac{1}{2}$, \eqref{dissipative} reduces to the well-known C\'{o}rdoba-C\'{o}rdoba-Fontelos model
\begin{eqnarray}\label{dissipative-CCF}
\partial_t\theta-H\theta\theta_x+\Lambda^\gamma\theta=0,
\end{eqnarray}
where $H\theta$ is the Hilbert transform of $\theta$. It was first proposed by C\'{o}rdoba, C\'{o}rdoba and Fontelos in \cite{[Cordoba-Cordoba-Fontelos05]} as a one-dimensional analogue of the two-dimensional surface quasigeostrophic equation (SQG) \cite{[Constantin-Majda-Tabak]}.
There have been a number of mathematical studies on the well-posedness for SQG.
We will not review here in detail the known results for the SQG and related equation. One can refer to \cite{[Caffarelli-Vasseur],[Chae-Constantin-Cordoba-Gancedo],[Constantin-Vicol],[Constantin-Wu],[Cordoba],[Cordoba-Martinez],[CotiZelati-Vicol],[Jeong-Kim],[Jeong-Kim-Yao],[Kiselev-Nazarov-Volberg],[Kiselev-Ryzhik-Yao-Zlato],[Kiselev-Yao-Zlatos],[Zlatos]} and the references therein  for more details and the recent progresses. Here we briefly summarize the progress related to the CCF model and related equations.
Concerning the inviscid CCF model, i.e., \eqref{M-CCF-T} with $n=1$ and $\alpha=\frac12$,
C\'{o}rdoba et al.\cite{[Cordoba-Cordoba-Fontelos05]} first obtained an ingenious nonlinear weighed inequality
 for the Hilbert transform by the use of Meillin transform and complex analysis, and proved
 smooth solutions must blow up in finite time for a generic family of even initial data (see also \cite{[Cordoba-Cordoba-Fontelos06]} for another proof of the blowup for the initial data not necessarily even).
This blow-up phenomena was later proved by Silvestre and Vicol in \cite{[Silvestre-Vicol]} via four essentially different methods.
Based on completely real-variable arguments, Li and Rodrigo \cite{[Li-Rodrigo20]} recently provided a short proof of the nonlinear inequality first proved by C\'{o}rdoba et al. in \cite{[Cordoba-Cordoba-Fontelos05]}, and obtained several new weighted inequalities for the Hilbert transform and various nonlinear versions, which can be applied to show the finite-time blow-up of smooth solutions to \eqref{M-CCF-T} with $n=1$ and $\alpha\in(0,1)$.
For the dissipative case, the authors in \cite{[Cordoba-Cordoba-Fontelos05]} also obtained the global well-posedness of \eqref{dissipative-CCF} for the positive $H^2$ initial data in the subcritical case $1<\gamma<2$. Later in \cite{[Dong2008]}, by adapting the method of continuity first used in \cite{[Kiselev-Nazarov-Volberg]}, Dong established the global regularity of solutions to \eqref{dissipative-CCF} in the critical case $\gamma=1$ for arbitrary initial data in appropriate critical Sobolev space. In \cite{[Li-Rodrigo08]}, Li and Rodrigo proved the finite-time blowup of smooth solutions to \eqref{dissipative-CCF} in the supercritical case $0<\gamma<\frac12$ (see also \cite{[Kiselev],[Li-Rodrigo20],[Silvestre-Vicol]} for different proofs for this blow-up result).
Recently, for each smooth nonnegative initial data, Ferreira and Moitinho \cite{[Ferreira-Moitinho20]} obtained the existence of global classical solutions to \eqref{dissipative-CCF} for $\gamma\in(\gamma_1,1)$ with $\gamma_1$ depending on the $H^{\frac32}$-norm of the initial data.
For the remaining case $\frac12\leq\gamma<1$, whether smooth solutions to \eqref{dissipative-CCF} may blow up in finite time is currently still open.
Some previously mentioned well-poseness results on the dissipative CCF model \eqref{dissipative-CCF} have been extended to the interpolation $\alpha$-CCF model, i.e., \eqref{dissipative} with $n=1$ and $0<\alpha<\frac12$, by Ferreira and Moitinho in a recent work \cite{[Ferreira-Moitinho]}.

A natural multi-dimensional generalization of the inviscid CCF equation, i.e., \eqref{M-CCF-T} with $n\geq2$ and $\alpha=\frac12$,  was considered in \cite{[Balodis-Cordoba]}, where
Balodis and C\'{o}rdoba presented a weighted nonlinear inequality for the Riesz transform by using Meillin transform and spherical harmonic expansion,
and  obtained the blow-up of smooth solution for any nonnegative, not-identically zero initial data. When $n=2$, such result was also proved for a similar equation in \cite{[Dong-Li]} independently. For the two-dimensional case with fractional dissipation, Li and Rodrigo \cite{[Li-Rodrigo09]} proved the finite-time blow-up of radial smooth  solutions to \eqref{dissipative} for $\frac{1}{4}\leq\alpha<1$ and $0<\gamma<\alpha$. Later, in \cite{[Dong]}, Dong was able to obtain the blow-up of smooth solutions to \eqref{M-CCF-T} with full range $\alpha\in (0,1)$ and $n\geq2$ for any smooth, radially symmetric and nonnegative initial data with compact support and its positive maximum attained at the origin.
Recently, motivated by \cite{[Silvestre-Vicol]}, Jiu and Zhang \cite{[Jiu-Zhang]} proved the finite time singularity of solutions to \eqref{M-CCF-T} with $0<\alpha<1$ and $n\geq2$ for smooth initial data $\theta_0$ with $\displaystyle\sup_{x\in\R^n}\theta_0(x)>0$ via the De Giorgi iteration technique.
This iteration strategy in \cite{[Silvestre-Vicol]} has also been adapted by Alonso-Or\'{a}n and Mart\'{i}nez \cite{[A-O-M]} to the proof of finite time blow-up for non-local active scalar equations \eqref{M-CCF-T} on compact Riemannian manifolds. Very recently, Li, Liu and Zhang \cite{[Li-Li-Zhang]} studied a related model, i.e., \eqref{M-CCF-T} with an additional power type of damping term:
\begin{equation}\label{M-CCF-T-damping}
\left\{\ba
&\partial_t\theta+u\cdot\nabla \theta +\kappa|\theta|^{\nu-1}\theta= 0, ~(x,t)\in \R^{n}\times\R_+,\\&u=\nabla\Lambda^{-2+2\alpha}\theta,\\
&\theta(x,0)=\theta_{0}(x),~x\in\R^n, \ea\ \right.
\end{equation}
where $n\geq1$, $0<\alpha<1$,   $\kappa\in \R$ and $\nu>0$. In \cite{[Li-Li-Zhang]}, by some change of time variable to implement the iteration technique in \cite{[Silvestre-Vicol]}, the authors showed that the damping term can not avoid the singularity formation in finite time and that for particular initial data depending usually on the size of their $L^1$ and $L^\infty$-norm, the solutions to \eqref{M-CCF-T-damping} must blow up, independent of the value of $\kappa$.

In this paper, we revisit the problem of finite-time blowup for \eqref{M-CCF-T} with $n\geq2$ and $0<\alpha<1$.
Our main result on \eqref{M-CCF-T} can be stated as
\begin{theorem}\label{the-1}
Let the initial data $\theta_0$ be a radial Schwartz function.
There exists a constant $A(n,\alpha)>0$ depending only on $n$ and $\alpha$ such that if
\begin{eqnarray}\label{Class-initial-data}
\int_{\R^n}\frac{\theta_0(0)-\theta_0(x)}{|x|^n}e^{-|x|}dx> A(n,\alpha)\|\theta_0\|_{L^\infty},
\end{eqnarray}
then the smooth solution $\theta$ to \eqref{M-CCF-T} blows up in finite time.
\end{theorem}
\begin{remark}
We may construct the radial Schwartz initial data $\theta_0$ fulfilling the inequality \eqref{Class-initial-data} as a smooth bump function below:
\begin{eqnarray*}
\theta_0(x)
=
\begin{cases}
1,
&  \mbox{if $|x|\leq \delta,$ }\\
\frac{e^{-\frac{1}{2\delta-|x|}}}{e^{-\frac{1}{2\delta-|x|}}+e^{-\frac{1}{|x|-\delta}}},
& \mbox{if $\delta<|x|<2\delta,$ } \\
0,
&  \mbox{if $|x|\geq 2\delta,$ }\\
\end{cases}
\end{eqnarray*}
where the constant $\delta\in\Big(0,\frac12 e^{-\frac{eA(n,\alpha)}{\omega_{n-1}}}\Big)$ is arbitrarily fixed. Then, by a integration by parts, we have
\begin{align*}
\int_{\mathbb{R}^n}\frac{\theta_0(0)-\theta_0(x)}{|x|^n}e^{-|x|}dx
&\geq\int_{|x|\geq2\delta }\frac{\theta_0(0)-\theta_0(x)}{|x|^n}e^{-|x|}dx=\omega_{n-1}\int^\infty_{2\delta}\frac{e^{-r}}{r}dr\\
&\geq\omega_{n-1}\int^1_{2\delta}\frac{e^{-r}}{r}dr=\omega_{n-1}\Big(\frac{\ln\frac{1}{2\delta}}{e^{2\delta}}
+\int^1_{2\delta}\frac{\ln r}{e^r}dr\Big)\\
&\geq\omega_{n-1}\Big(\frac{\ln\frac{1}{2\delta}}{e^{2\delta}}
+\ln 2\delta\int^1_{2\delta}\frac{dr}{e^r}\Big)\\
&=\frac{\omega_{n-1}}{e}\ln\frac{1}{2\delta}>A(n,\alpha)\|\theta_0\|_{L^\infty},
\end{align*}
which is the desired inequality \eqref{Class-initial-data}.
Here $\omega_{n-1}$ denotes the surface area of the uint sphere in $\mathbb{R}^n$.
\end{remark}
For the proof of Theorem \ref{the-1}, we consider the evolution of a weighted integral of the solution, and show that such quantity satisfies some ordinary differential inequality and will blow up in finite time (see \eqref{weighted-integral-solution} and \eqref{ordinary-differential-inequality}). We then obtain the blow-up of the solution. This type of strategy can be, for instance, found in \cite{[Cordoba-Cordoba-Fontelos05],[Dong],[Li-Rodrigo08],[Li-Rodrigo09],[Li-Rodrigo20]}. The key ingredient to adapt to this approach is to derive a weighted nonlinear inequality (see \eqref{nonlinear-inequality-exponential-weight}). Inspired by \cite{[Li-Rodrigo20]}, we prove the weighted nonlinear inequality \eqref{nonlinear-inequality-exponential-weight} by some completely real variable based techniques.
It will be seen later that \eqref{nonlinear-inequality-exponential-weight} is reduced to the estimate of the integral $\mathcal{I}$ (see \eqref{representation-2}).
For the case of $0<\alpha<\frac{1}{2}$, by discarding all the positive term in the integral, a lower bound of $\mathcal{I}$ can be readily obtained with help of the fact that the series $\displaystyle\sum^\infty_{k=0}a_{2k}(\alpha)$ is convergent in this case (see Corollary \ref{Positivity-Convergence}). The remaining case $\alpha\in[\frac12,1)$ is more involved and needs more refined arguments (see \eqref{3.13}). The difficulty in both cases is that the series $\displaystyle\sum^\infty_{k=0}a_{2k}(\alpha)$ diverges for $\alpha\in[\frac12,1)$. By some detailed singular integral estimates and the fact that the series $\displaystyle\sum^\infty_{k=0}\frac{a_{2k}(\alpha)}{k}$ is convergent for the full range $\alpha\in(0,1)$, we get the desired weighted inequality.
In comparison with \cite{[Dong]}, a different approach of the proof of the finite-time blow-up of smooth solution to \eqref{M-CCF-T} is given in this paper.

Throughout this paper, we will use $C$ to denote a positive constant, whose value may change from line to line, and write $C_{n,\alpha}$ or $C(n,\alpha)$ to emphasize the dependence of a constant on $n$ and $\alpha$.
For $p\in[1,\infty]$, we denote $L^p(\R^n)$ the standard $L^p$-space and its  norm by  $\|\cdot\|_{L^p(\R^n)}$.
For $s\geq0$, we use the notation $H^s(\R^n)$ to denote the nonhomogeneous Sobolev space of $s$ order,
whose endowed norm is denoted by $\|\cdot\|_{H^s(\R^n)}=\|\cdot\|_{L^2(\R^n)}+\|\Lambda^s(\cdot)\|_{L^2(\R^n)}$ (see \cite{[Bahouri-Chemin-Danchi]} for more details).
For a sake of the convenience, the $L^{p}(\R^{n})$-norm of a function $f$ is always abbreviated as $\|f\|_{L^p}$ and
the
$H^s(\R^n)$-norm as
$\norm{f}_{H^s}$.
The real Gamma function $\Gamma$ is defined by
\begin{eqnarray*}
\Gamma(s)=\int^\infty_0t^{s-1}e^{-t}dt,~{\rm for}~s>0.
\end{eqnarray*}
The Gamma function can be expressed as a limit of some sequence \cite{[Rudin]}, that is,
for $s>0$,
\begin{eqnarray}\label{Gamma-limit}
\Gamma(s)
=\lim_{k\rightarrow+\infty}\frac{k^sk!}{s(s+1)\cdot\cdot\cdot(s+k)}.
\end{eqnarray}
The related real Beta function $B$ is defined by
\begin{eqnarray*}
B(p,q)=\int_0^1t^{p-1}(1-t)^{q-1}dt,~{\rm for}~p>0,~q>0.
\end{eqnarray*}
It is well-known that
\begin{eqnarray*}
B(p,q)=\frac{\Gamma(p)\Gamma(q)}{\Gamma(p+q)},~ p>0,~q>0.
\end{eqnarray*}
Finally, let $\mathbb{S}^{n-1}$ be the unit sphere in $\mathbb{R}^n$, i.e., $\mathbb{S}^{n-1}=\{x\in\mathbb{R}^n:|x|=1\}$
and $\omega_{n-1}$ be its surface area. We recall that
\begin{eqnarray*}
\omega_{n-1}=\frac{2\pi^{\frac{n}{2}}}{\Gamma(\frac{n}{2})}.
\end{eqnarray*}

The remaining part of this paper is organized as follows. In Section 2, we first recall the local well-posedness for the model \eqref{M-CCF-T}, and then present some useful facts utilized later in this paper.
Section 3 is devoted to two weighted nonlinear inequalities related to the model \eqref{M-CCF-T}.
The proof of Theorems \ref{the-1} is given in Section 4.
\section{Preliminaries}
In this section, we will present some basic useful facts needed later.
We begin with the local well-posedness of \eqref{M-CCF-T} in the Sobolev space $H^s$ for some appropriate $s>0$ , whose proof was given in \cite{[Chae]}.
\begin{lemma}\label{local-well-posedness}
{\rm(1)} Let $n\geq2$ and $\frac{1}{2}<\alpha<1$. Then for each $\theta_0\in H^s$ with $s>\frac{n}{2}+2$, there exists a positive $T=T(\|\theta_0\|_{H^s})$ such that the equation \eqref{M-CCF-T} admits a unique solution $\theta$ in $C([0,T); H^s)\cap {{\rm Lip}}((0,T);H^{s-1}).$

{\rm(2)} Let $n\geq2$ and $0<\alpha\leq\frac{1}{2}$. Then for each $\theta_0\in H^s$ with $s>\frac{n}{2}+1$, there exists a positive $T=T(\|\theta_0\|_{H^s})$ such that the equation \eqref{M-CCF-T} admits a unique solution $\theta$ in $C([0,T); H^s)\cap {\rm Lip}((0,T);H^{s-1})$. Furthermore, if $T^\ast$ is the first time the solution cannot be continued in $C([0,T^\ast); H^s)$, then there necessarily holds $$\int_0^{T^\ast}\|(\mathcal{R}\otimes\mathcal{R})\Lambda^{2\alpha}\theta(\cdot,t)\|_{L^\infty}dt=\infty,$$
where $\Big((\mathcal{R}\otimes\mathcal{R})(f)\Big)_{jk}=\mathcal{R}_j\mathcal{R}_k(f)$ denotes the tensor product of the Riesz transform.
\end{lemma}
The next lemma shows that the radial symmetry of the initial data is preserved by the solution to \eqref{M-CCF-T}.
\begin{lemma}\label{radial-property-preserved}
If $\theta$ is a smooth solution to \eqref{M-CCF-T} with the radially symmetric initial data $\theta_0$, then  $\theta(x,t)$ is radially symmetric for all $t>0$ in its life span.
\end{lemma}
\textbf{Proof}. Let $\mathbf{O}\in\mathbb{R}^{n\times n}$ be any orthogonal matrix. By the uniqueness of solutions to \eqref{M-CCF-T} and  the radial property of $\theta_0$, it suffices to show that the function
\begin{eqnarray*}
\theta_{\mathbf{O}}(x,t):=\theta(\mathbf{O}x,t)
\end{eqnarray*}
is also a solution to \eqref{M-CCF-T} with the initial data $\theta_0(\mathbf{O}x)$. Indeed, standard computations give that
\begin{eqnarray*}
(\partial_t\theta_{\mathbf{O}})(x,t)=\partial_t(\theta(\mathbf{O}x,t))=(\partial_t\theta)(\mathbf{O}x,t)
\end{eqnarray*}
and
\begin{eqnarray*}
(\nabla_x\theta_{\mathbf{O}})(x,t)=\nabla_x(\theta(\mathbf{O}x,t))={\mathbf{O}}^T(\nabla_x\theta)(\mathbf{O}x,t).
\end{eqnarray*}
By the second equation in \eqref{M-CCF-T}, the integral representation of the Riesz potential and integration by parts, we can derive that
\begin{align*}
u_{\mathbf{O}}(x,t)
&=C_{n,\alpha}P.V.\int_{\R^{n}}\frac{x-y}{|x-y|^{n+2\alpha}}\theta(\mathbf{O}y, t)dy\\
&=C_{n,\alpha}P.V.\int_{\R^{n}}\frac{x-{\mathbf{O}}^{-1}z}{|x-{\mathbf{O}}^{-1}z|^{n+2\alpha}}\theta(z, t)dz\\
&={\mathbf{O}}^{-1}C_{n,\alpha}P.V.\int_{\R^{n}}\frac{\mathbf{O} x-z}{|\mathbf{O} x-z|^{n+2\alpha}}\theta(z,t)dz\\
&={\mathbf{O}}^{-1} u(\mathbf{O} x,t),
\end{align*}
where
\begin{eqnarray*}
C_{n,\alpha}
=-\frac{2^{2\alpha-1}\Gamma(\frac{n}{2}+\alpha)}{\pi^{\frac{n}{2}}\Gamma(1-\alpha)}.
\end{eqnarray*}
Thus, we obtain
\begin{align*}
(\partial_t\theta_{\mathbf{O}}+u_{\mathbf{O}}\cdot\nabla\theta_{\mathbf{O}})(x,t)
&=(\partial_t\theta)(\mathbf{O}x,t)+{\mathbf{O}}^{-1} u(\mathbf{O} x,t)\cdot {\mathbf{O}}^T(\nabla\theta)(\mathbf{O}x,t)\\
&=(\partial_t\theta)(\mathbf{O}x,t)+u(\mathbf{O} x,t)\cdot \mathbf{O}{\mathbf{O}}^T(\nabla\theta)(\mathbf{O}x,t)\\
&=(\partial_t\theta+u\cdot\nabla\theta)(\mathbf{O}x,t)=0.
\end{align*}
The proof of Lemma \ref{radial-property-preserved} is then finished. \hfill\hfill$\square$\vskip12pt

We proceed to present a string of simple inequalities needed later.
\begin{lemma}\label{Young-inequality}
Let $\frac12<\alpha<1$. Given $f:[0,\infty)\mapsto\R$ a smooth and bounded function. We have
\begin{eqnarray*}
\int^\infty_0\frac{(f(\rho)-f(0))^2}{\rho^{2\alpha}}d\rho
\leq\epsilon\int^\infty_0\frac{(f(\rho)-f(0))^2}{\rho^{1+2\alpha}}d\rho
+
\Big(\frac{1}{(2-2\alpha)\epsilon}+\frac{4}{2\alpha-1}\Big)\|f\|^2_{L^\infty},
\end{eqnarray*}
\begin{eqnarray*}
\int_0^\infty\frac{(f(\rho)-f(0))^2}{\rho e^{\frac\rho2}}d\rho
\leq\epsilon\int_0^\infty\frac{(f(\rho)-f(0))^2}{\rho^{2}}d\rho+\Big(8+\frac1\epsilon\Big)\|f\|^2_{L^\infty}
\end{eqnarray*}
and
\begin{eqnarray*}
\int_0^\infty\frac{(f(\rho)-f(0))^2}{\rho^{2}}(1-e^{-\frac{\rho}{2}})d\rho
\leq
\epsilon\int_0^\infty\frac{(f(\rho)-f(0))^2}{\rho^{2}}d\rho+\Big(4+\frac{1}{4\epsilon}\Big)\|f\|^2_{L^\infty},
\end{eqnarray*}
for any positive $\epsilon$.
\end{lemma}
\textbf{Proof.} Since $\frac{1}{2}<\alpha<1$, by H\"{o}lder's and Young's inequality, we first derive that
\begin{align*}
\int\limits^\infty_0\frac{(f(\rho)-f(0))^2}{\rho^{2\alpha}}d\rho
&\leq2\|f\|_{L^\infty}\Big[\int\limits^1_0\frac{(f(\rho)-f(0))^2}{\rho^{1+2\alpha}}d\rho\Big]^{\frac12}
\Big[\int^1_0\frac{d\rho}{\rho^{2\alpha-1}}\Big]^{\frac12}
+4\|f\|^2_{L^\infty}\int\limits^\infty_1\frac{d\rho}{\rho^{2\alpha}}\\
&\leq\epsilon\int^\infty_0\frac{(f(\rho)-f(0))^2}{\rho^{1+2\alpha}}d\rho
+\Big(\frac{1}{(2-2\alpha)\epsilon}+\frac{4}{2\alpha-1}\Big)\|f\|^2_{L^\infty},
\end{align*}
where $\epsilon>0$ is any positive constant.

Similarly, also by the elementary inequalities
\begin{eqnarray*}
e^{\frac{\rho}{2}}>\frac{\rho}{2}~{\rm and}~1-e^{-\frac{\rho}{2}}<\frac{\rho}{2},~ {\rm for}~\rho>0,
\end{eqnarray*}
we can obtain that
\begin{align*}
\int_0^\infty\frac{(f(\rho)-f(0))^2}{\rho e^{\frac\rho2}}d\rho
&\leq2\|f\|_{L^\infty}\int_0^1\frac{|f(\rho)-f(0)|}{\rho}d\rho
+8\|f\|^2_{L^\infty}\int_1^\infty\frac{d\rho}{\rho^2}\\
&\leq2\|f\|_{L^\infty}\Big(\int_0^1\frac{|f(\rho)-f(0)|^2}{\rho^2}d\rho\Big)^{\frac12}
+8\|f\|^2_{L^\infty}\\
&\leq\epsilon\int_0^\infty\frac{(f(\rho)-f(0))^2}{\rho^{2}}d\rho+\Big(8+\frac1\epsilon\Big)\|f\|^2_{L^\infty}.
\end{align*}
and
\begin{align*}
\int\limits_0^\infty\frac{(f(\rho)-f(0))^2}{\rho^{2}}(1-e^{-\frac{\rho}{2}})d\rho
&\leq\|f\|_{L^\infty}\int\limits_0^1\frac{|f(\rho)-f(0)|}{\rho}d\rho
+4\|f\|^2_{L^\infty}\int\limits_1^\infty\frac{d\rho}{\rho^{2}}\\
&\leq\|f\|_{L^\infty}\Big(\int_0^1\frac{(f(\rho)-f(0))^2}{\rho^{2}}d\rho\Big)^{\frac12}
+4\|f\|^2_{L^\infty}\\
&\leq\epsilon\int_0^\infty\frac{(f(\rho)-f(0))^2}{\rho^{2}}d\rho+\Big(4+\frac{1}{4\epsilon}\Big)\|f\|^2_{L^\infty}.
\end{align*}
We then finish the proof of Lemma \ref{Young-inequality}. \hfill\hfill$\square$\vskip12pt

Finally, we end with this section by introducing an auxiliary function defined by
\begin{eqnarray}\label{g-definition}
g(\lambda,m,\beta)\triangleq\int_0^\pi\frac{\sin^m\mu d\mu}{(1-2\lambda\cos\mu+\lambda^2)^{\frac{m}{2}+\beta}},
\end{eqnarray}
where $\lambda\in[0,\infty)$, $m\in\mathbb{N}$ and $\beta>0$.
The motivation for this function will be clear later in the proof of Proposition \ref{n-nonlinear-inequality}. One can extend $g$ to $\lambda\in(-\infty,\infty)$ by using \eqref{g-definition} and readily see that $g$ is even, infinitely differentiable (and real analytic) in $\lambda\in(-1,1)$ and its Taylor's series at $0$ has the form:
\begin{eqnarray*}
g(\lambda,m,\beta)=\sum^\infty_{k=0}a_{2k}(m,\beta)\lambda^{2k}.
\end{eqnarray*}
The following lemma is mainly from \cite{[Dong]}.
Since it is vital for our arguments and is also of independent interest, we include the proof here for the sake of completeness.
\begin{lemma}\label{Taylor-series}
For $k\geq0$, it holds that $a_{2k}(m,\beta)>0$ and
\begin{eqnarray*}
\displaystyle\lim_{k\rightarrow\infty}\frac{a_{2k}(m,\beta)}{k^{2\beta-2}}
=\frac{\Gamma(\frac12)\Gamma(\frac{m}{2}+\frac12)}{\Gamma(\beta)\Gamma(\frac{m}{2}+\beta)}.
\end{eqnarray*}
\end{lemma}
\textbf{Proof.} Direct differentiation yields that
\begin{eqnarray*}
\frac{\partial g}{\partial \lambda}(\lambda,m,\beta)
=-(m+2\beta)\int_0^\pi\frac{\sin^m\mu(\lambda-\cos\mu)}{(1-2\lambda\cos\mu+\lambda^2)^{\frac{m}{2}+\beta+1}}d\mu
\end{eqnarray*}
and
\begin{align*}
\frac{\partial^2 g}{\partial \lambda^2}(\lambda,m,\beta)
&=-(m+2\beta)\int_0^\pi\frac{\sin^m\mu d\mu}{(1-2\lambda\cos\mu+\lambda^2)^{\frac{m}{2}+\beta+1}}\\
&\ \ \
+(m+2\beta)(m+2\beta+2)\int_0^\pi\frac{\sin^m\mu(\lambda-\cos\mu)^2 }{(1-2\lambda\cos\mu+\lambda^2)^{\frac{m}{2}+\beta+2}}d\mu.
\end{align*}
Noting that
\begin{eqnarray*}
(\lambda-\cos\mu)^2 =(1-2\lambda\cos\mu+\lambda^2)-\sin^2\mu,
\end{eqnarray*}
we further derive that
\begin{eqnarray*}
\frac{\partial^2 g}{\partial \lambda^2}(\lambda,m,\beta)
=(m+2\beta)\Big[(m+2\beta+1)g(\lambda,m,\beta+1)
 -(m+2\beta+2)g(\lambda,m+2,\beta+1)\Big],
\end{eqnarray*}
which implies that
\begin{eqnarray}\label{recurrence}
\begin{split}
\frac{\partial^{2k} g}{\partial \lambda^{2k}}(0,m,\beta)
 &=(m+2\beta)\Big[(m+2\beta+1)\frac{\partial^{2k-2} g}{\partial \lambda^{2k-2}}(0,m,\beta+1)\\
 &\ \ \ \ \ \ \ \ \ \ \ \ \ \ \ \ \ \
 -(m+2\beta+2)\frac{\partial^{2k-2} g}{\partial \lambda^{2k-2}}(0,m+2,\beta+1)\Big],
\end{split}
\end{eqnarray}
for any integer $k\geq1$.
Since
\begin{eqnarray*}
a_0(m,\beta)=g(0,m,\beta)=2\int_0^{\frac{\pi}{2}}\sin^m\mu d\mu=B\Big(\frac12,\frac{m+1}{2}\Big)>0,
\end{eqnarray*}
then, from the recurrence relation \eqref{recurrence}, by induction, we can obtain that, for $k\geq1$,
\begin{align*}
\frac{\partial^{2k} g}{\partial \lambda^{2k}}(0,m,\beta)
 &=B\Big(\frac12,\frac{m+1}{2}\Big)(2k-1)!!\cdot2\beta(2\beta+2)\cdot\cdot\cdot(2\beta+2k-2)\\
 &\ \ \ \
 \cdot\frac{(m+2\beta)(m+2\beta+2)\cdot\cdot\cdot(m+2\beta+2k-2)}{(m+2)(m+4)\cdot\cdot\cdot(m+2k)}.
\end{align*}
It follows from
\begin{eqnarray*}
a_{2k}(m,\beta)=\frac{1}{(2k)!}\frac{\partial^{2k} g}{\partial \lambda^{2k}}(0,m,\beta)
\end{eqnarray*}
that, for $k\geq1$,
\begin{eqnarray*}
a_{2k}(m,\beta)
=B\Big(\frac12,\frac{m+1}{2}\Big)
\frac{(\beta)_k(\frac{m}{2}+\beta)_k}{k!(\frac{m}{2}+1)_k}>0,
\end{eqnarray*}
where the Pochhammer's symbol's $(x)_k$ is defined by
\begin{eqnarray*}
(x)_k
=
\begin{cases}
1,
&  \mbox{if $k=0,$ }\\
x(x+1)\cdot\cdot\cdot(x+k-1),
& \mbox{if $k\geq1.$ } \\
\end{cases}
\end{eqnarray*}
It follows from \eqref{Gamma-limit} that
\begin{align*}
\displaystyle\lim_{k\rightarrow\infty}\frac{a_{2k}(m,\beta)}{k^{2\beta-2}}
&=
B\Big(\frac12,\frac{m+1}{2}\Big)
\displaystyle\lim_{k\rightarrow\infty}\Big[\frac{\frac{k^{\frac{m}{2}+1}k!}{(\frac{m}{2}+1)_{k+1}}}{\frac{k^\beta k!}{(\beta)_{k+1}}\cdot\frac{k^{\frac{m}{2}+\beta}k!}{(\frac{m}{2}+\beta)_{k+1}}}
\cdot
\frac{k(\frac{m}{2}+1+k)}{(\beta+k)(\frac{m}{2}+\beta+k)}\Big]\\
&=B\Big(\frac12,\frac{m+1}{2}\Big)\frac{\Gamma(\frac{m}{2}+1)}{\Gamma(\beta)\Gamma(\frac{m}{2}+\beta)}
=\frac{\Gamma(\frac12)\Gamma(\frac{m}{2}+\frac12)}{\Gamma(\beta)\Gamma(\frac{m}{2}+\beta)}.
\end{align*}
The proof of Lemma \ref{Taylor-series} is then complete. \hfill\hfill$\square$\vskip12pt
For the convenience of our later application, we denote
\begin{eqnarray}\label{G-definition}
G_\alpha(\lambda)
\triangleq\int_0^\pi\frac{\sin^n\mu d\mu}{(1-2\lambda\cos\mu+\lambda^2)^{\frac{n}{2}+\alpha}}
=\sum^\infty_{k=0}a_{2k}(\alpha)\lambda^{2k},
\end{eqnarray}
where $\lambda\in(-1,1)$, $n\geq2$ and $\alpha\in(0,1)$.
Then, for any $\alpha\in(0,1)$,
\begin{eqnarray}\label{G-0-value}
G_\alpha(0)=B\Big(\frac12,\frac{n+1}{2}\Big)
\end{eqnarray}
and
it is easy to check that, for any $\lambda>0$ and $\lambda\neq1$ (notice that $G_\alpha(\lambda)$ is not defined at 1 for $\alpha\in[\frac12,1)$),
\begin{eqnarray}\label{transform-identity}
G_\alpha\Big(\frac{1}{\lambda}\Big)=\lambda^{n+2\alpha}G_\alpha(\lambda).
\end{eqnarray}
In addition, $G_\alpha$ is smooth in $[0,1)\cup(1,\infty)$ with a possible singularity at $1$. The singularity is of order $|\lambda-1|^{1-2\alpha}$
when $\alpha\in(\frac12,1)$ and of order $\log|\lambda-1|$ when $\alpha=\frac12$.

As a corollary of Lemma \ref{Taylor-series}, we immediately have
\begin{coro}\label{Positivity-Convergence}
For $n\geq2$ and $\alpha\in(0,1)$, $G_\alpha(\lambda)$, $G'_\alpha(\lambda)$ and $G''_\alpha(\lambda)$ are all positive in $\lambda\in(0,1)$.
 Furthermore, the series $\displaystyle\sum^\infty_{k=0}a_{2k}(\alpha)$ is convergent if and only if $\alpha\in(0,\frac12)$, and $\displaystyle\sum^\infty_{k=1}\frac{a_{2k}(\alpha)}{k}$ is convergent for $\alpha\in(0,1)$.
\end{coro}
\textbf{Proof.} By Lemma \ref{Taylor-series}, we have, for $0<\alpha<1$,
\begin{eqnarray*}
\displaystyle\lim_{k\rightarrow\infty}\frac{a_{2k}(\alpha)}{k^{2\alpha-2}}
=\frac{\Gamma(\frac12)\Gamma(\frac{n}{2}+\frac12)}{\Gamma(\alpha)\Gamma(\frac{n}{2}+\alpha)},
\end{eqnarray*}
which implies that $\displaystyle\sum^\infty_{k=0}a_{2k}(\alpha)$ is convergent if and only if $0<\alpha<\frac12$, and that $\displaystyle\sum^\infty_{k=1}\frac{a_{2k}(\alpha)}{k}$ is convergent for $\frac12\leq\alpha<1$. In addition, it follows that the radius of the convergence of \eqref{G-definition} is 1 for $0<\alpha<1$. Thus, by differentiation term by term, we ge that, for $\lambda\in(-1,1)$,
\begin{eqnarray}\label{1-order}
G'_\alpha(\lambda)=\sum^\infty_{k=1}2ka_{2k}(\alpha)\lambda^{2k-1}
\end{eqnarray}
and
\begin{eqnarray*}
G''_\alpha(\lambda)=\sum^\infty_{k=1}2k(2k-1)a_{2k}(\alpha)\lambda^{2k-2},
\end{eqnarray*}
which along with $a_{2k}(\alpha)>0$ and \eqref{G-definition} that $G_\alpha(\lambda)$, $G'_\alpha(\lambda)$ and $G''_\alpha(\lambda)$ are all positive in $\lambda\in(0,1)$. The proof of Corollary \ref{Positivity-Convergence} is then finished. \hfill\hfill$\square$\vskip12pt
\section{Two weighted inequalities for nonlinear term}
In this section, we will prove two nonlinear weighted inequalities for the model \eqref{M-CCF-T}. The one-dimensional analogous inequalities were established in \cite{[Li-Rodrigo20]}.
By abuse of notation, we will not distinguish $f(x)$ and $f(|x|)$ for a radially symmetric function $f$.
\begin{proposition}\label{n-nonlinear-inequality}
Let $n\geq2$ and $0<\alpha<1$. Let $f:\R^n\rightarrow\R$ be a radial Schwartz function. Then
\begin{eqnarray*}
\int_{\R^n}\frac{\Lambda^{-2+2\alpha}\nabla f(x)\cdot\nabla f(x)}{|x|^{n}}dx
\geq
\frac{\alpha2^{2\alpha-1}\Gamma(\frac{n}{2}+\alpha)}{\Gamma(1-\alpha)\Gamma(\frac{n}{2}+1)}
\int_{\R^n}\frac{(f(0)-f(x))^2}{|x|^{n+2\alpha}}dx.
\end{eqnarray*}
\end{proposition}

\textbf{Proof.} By the integral representation of the Riesz potential \cite{[Stein]} and the radial assumption on $f$, we have
\begin{align*}
\Lambda^{-2+2\alpha}(\nabla f)(x)
&=c_{n,\alpha}\int_{\R^{n}}\frac{\nabla f(y)}{|x-y|^{n-2+2\alpha}}dy\\
&=c_{n,\alpha}\int_{\R^{n}}\frac{f'(|y|)}{|x-y|^{n-2+2\alpha}}\frac{y}{|y|}dy\\
&=c_{n,\alpha}\int_0^\infty f'(\rho)\rho^{n-1}\Big[\int_{\mathbb{S}^{n-1}}\frac{zd\sigma(z)}{|x-\rho z|^{n-2+2\alpha}}\Big]d\rho,
\end{align*}
where
\begin{eqnarray*}
c_{n,\alpha}=\frac{\Gamma(\frac{n}{2}-1+\alpha)}{\pi^{\frac{n}{2}}2^{2-2\alpha}\Gamma(1-\alpha)}.
\end{eqnarray*}
For any $x\in\mathbb{R}^n\backslash\{0\}$, pick an orthogonal matrix $\mathbf{O}\in\mathbb{R}^{n\times n}$ such that $\mathbf{O}e_1=\frac{x}{|x|}$, where $e_1=(1,0,...,0)$.
It follows that, for any $x\in\mathbb{R}^n\backslash\{0\}$,
\begin{align*}
\Lambda^{-2+2\alpha}\nabla f(x)\cdot\nabla f(x)
&=c_{n,\alpha}f'(|x|)\int_0^\infty f'(\rho)\rho^{n-1}\Big(\int_{\mathbb{S}^{n-1}}\frac{z\cdot\frac{x}{|x|}d\sigma(z)}{|x-\rho z|^{n-2+2\alpha}}\Big)d\rho\\
&=c_{n,\alpha}f'(|x|)\int_0^\infty f'(\rho)\rho^{n-1}\Big(\int_{\mathbb{S}^{n-1}}\frac{z\cdot \mathbf{O}e_1d\sigma(z)}{||x|\mathbf{O}e_1-\rho z|^{n-2+2\alpha}}\Big)d\rho\\
&=c_{n,\alpha}f'(|x|)\int_0^\infty f'(\rho)\rho^{n-1}\Big(\int_{\mathbb{S}^{n-1}}\frac{e_1\cdot {\mathbf{O}}^{-1}zd\sigma(z)}{||x|e_1-\rho {\mathbf{O}}^{-1}z|^{n-2+2\alpha}}\Big)d\rho\\
&=c_{n,\alpha}f'(|x|)\int_0^\infty f'(\rho)\rho^{n-1}\Big(\int_{\mathbb{S}^{n-1}}\frac{z_1d\sigma(z)}{||x|e_1-\rho z|^{n-2+2\alpha}}\Big)d\rho.
\end{align*}
By a change of variables formula (see e.g., pp. 592 of \cite{[Grafakos]}) and integration by parts, we obtain that, for $\rho\neq|x|$,
\begin{align*}
\int_{\mathbb{S}^{n-1}}\frac{z_1d\sigma(z)}{||x|e_1-\rho z|^{n-2+2\alpha}}
&=\int_{-1}^1\int_{\sqrt{1-s^2}\mathbb{S}^{n-2}}\frac{sd\sigma(z)}{((|x|-\rho s)^2+\rho^2|z|^2)^{\frac{n}{2}-1+\alpha}}
\frac{ds}{\sqrt{1-s^2}}\\
&=\omega_{n-2}\int_{-1}^1\frac{s(1-s^2)^{\frac{n-3}{2}}ds}{(|x|^2-2|x|\rho s+\rho^2)^{\frac{n}{2}-1+\alpha}}\\
&=\omega_{n-2}\int_0^\pi\frac{\cos\mu\sin^{n-2}\mu d\mu}{(|x|^2-2|x|\rho\cos\mu+\rho^2)^{\frac{n}{2}-1+\alpha}}\\
&=\frac{(n-2+2\alpha)\omega_{n-2}}{n-1}\rho|x|\int_0^\pi\frac{\sin^n\mu d\mu}{(|x|^2-2|x|\rho\cos\mu+\rho^2)^{\frac{n}{2}+\alpha}}\\
&=\frac{(n-2+2\alpha)\omega_{n-2}\rho}{(n-1)|x|^{n-1+2\alpha}}
\int_0^\pi\frac{\sin^n\mu d\mu}{(1-2\frac{\rho}{|x|}\cos\mu+\frac{\rho^2}{|x|^2})^{\frac{n}{2}+\alpha}}\\
&=\frac{(n-2+2\alpha)\omega_{n-2}}{n-1}
\frac{\rho}{|x|^{n-1+2\alpha}}
G_\alpha\Big(\frac{\rho}{|x|}\Big),
\end{align*}
where $\omega_{n-2}=\frac{2\pi^{\frac{n-1}{2}}}{\Gamma(\frac{n-1}{2})}$ is the surface area of $\mathbb{S}^{n-2}$
and the function $G_\alpha$ is defined as \eqref{G-definition}.
Therefore, integrating by parts along with the boundary conditions
\begin{eqnarray*}
\lim_{\rho\rightarrow0^+}
\rho^nG_\alpha\Big(\frac{\rho}{|x|}\Big)(f(\rho)-f(|x|))=0,
\end{eqnarray*}
\begin{eqnarray*}
\lim_{\rho\rightarrow+\infty}\rho^nG_\alpha\Big(\frac{\rho}{|x|}\Big)(f(\rho)-f(|x|))
=|x|^{n+2\alpha}\lim_{\rho\rightarrow+\infty}G_\alpha\Big(\frac{|x|}{\rho}\Big)\frac{f(\rho)-f(|x|)}{\rho^{2\alpha}}
=0,
\end{eqnarray*}
we have that, for any $x\in\mathbb{R}^n\backslash\{0\}$,
\begin{eqnarray}\label{expression-nonlinear-term}
\begin{split}
\Lambda^{-2+2\alpha}\nabla f(x)\cdot\nabla f(x)
&=c'_{n,\alpha}\frac{f'(|x|)}{|x|^{n-1+2\alpha}}
\int_0^\infty f'(\rho)\rho^nG_\alpha\Big(\frac{\rho}{|x|}\Big)d\rho\\
&=c'_{n,\alpha}\frac{f'(|x|)}{|x|^{n-1+2\alpha}}
\int_0^\infty \rho^nG_\alpha\Big(\frac{\rho}{|x|}\Big)\frac{d}{d\rho}(f(\rho)-f(|x|))d\rho\\
&=-c'_{n,\alpha}\frac{f'(|x|)}{|x|^{n-1+2\alpha}}
\int_0^\infty \frac{\partial}{\partial\rho}\Big(\rho^nG_\alpha\Big(\frac{\rho}{|x|}\Big)\Big)(f(\rho)-f(|x|))d\rho,
\end{split}
\end{eqnarray}
with
\begin{eqnarray*}
c'_{n,\alpha}
=c_{n,\alpha}\frac{(n-2+2\alpha)\omega_{n-2}}{n-1}
=\frac{2^{2\alpha-1}\Gamma(\frac{n}{2}+\alpha)}{\pi^{\frac12}\Gamma(1-\alpha)\Gamma(\frac{n+1}{2})}.
\end{eqnarray*}
Here, for the case of $\alpha\in[\frac12,1)$, we used the singularity order of $G_\alpha$ at $1$ to
show that when conducting integration by parts, the boundary term at the singularity $|x|$ vanishes.

It then follows from Fubini's theorem and integration by parts that
\begin{eqnarray}\label{representation-1}
\begin{split}
&\int_{\R^n}\frac{\Lambda^{-2+2\alpha}\nabla f(x)\cdot\nabla f(x)}{|x|^{n}}dx\\
&=c''_{n,\alpha}\int_0^\infty\frac{f'(r)}{r^{n+2\alpha}}
\Big(-\int_0^\infty \frac{\partial}{\partial\rho}\Big(\rho^nG_\alpha\Big(\frac{\rho}{r}\Big)\Big)(f(\rho)-f(r))d\rho\Big)dr\\
&=\frac{c''_{n,\alpha}}{2}\int_0^\infty\int_0^\infty\frac{1}{r^{n+2\alpha}}
 \frac{\partial}{\partial\rho}\Big(\rho^nG_\alpha\Big(\frac{\rho}{r}\Big)\Big)\frac{\partial}{\partial r}\Big((f(\rho)-f(r))^2\Big)drd\rho\\
&=\frac{c''_{n,\alpha}}{2}\int_0^\infty\Big(\frac{1}{r^{n+2\alpha}}\frac{\partial}{\partial\rho}\Big(\rho^nG_\alpha\Big(\frac{\rho}{r}\Big)\Big)
(f(\rho)-f(r))^2\Big)\Big|^{\infty}_{r\rightarrow0^+}d\rho\\
&\ \ \
-\frac{c''_{n,\alpha}}{2}
\int_0^\infty\int_0^\infty
\frac{\partial}{\partial r}\Big(\frac{1}{r^{n+2\alpha}}\frac{\partial}{\partial\rho}\Big(\rho^nG_\alpha\Big(\frac{\rho}{r}\Big)\Big)\Big)(f(\rho)-f(r))^2d\rho dr
\end{split}
\end{eqnarray}
with
\begin{eqnarray}\label{constant-11}
c''_{n,\alpha}=c'_{n,\alpha}\omega_{n-1}.
\end{eqnarray}
Here, for the case of $\alpha\in[\frac12,1)$, the singularity order of $G_\alpha$ at $1$ also ensures that the boundary term at the singularity $\rho$ is zero.

By \eqref{G-0-value}, we note that
\begin{eqnarray}\label{boundary-1}
\begin{split}
&\lim_{r\rightarrow\infty}\frac{(f(\rho)-f(r))^2}{r^{n+2\alpha}}\frac{\partial}{\partial\rho}\Big(\rho^nG_\alpha\Big(\frac{\rho}{r}\Big)\Big)\\
&=\lim_{r\rightarrow\infty}\frac{(f(\rho)-f(r))^2}{r^{n+2\alpha}}\Big[n\rho^{n-1}G_\alpha\Big(\frac{\rho}{r}\Big)+
\frac{\rho^n}{r}G'_\alpha\Big(\frac{\rho}{r}\Big)\Big]=0
\end{split}
\end{eqnarray}
and
\begin{eqnarray}\label{condition}
\begin{split}
&\lim_{r\rightarrow0^+}\frac{(f(\rho)-f(r))^2}{r^{n+2\alpha}}\frac{\partial}{\partial\rho}\Big(\rho^nG_\alpha\Big(\frac{\rho}{r}\Big)\Big)\\
&=(f(\rho)-f(0))^2\lim_{r\rightarrow0^+}\frac{\partial}{\partial\rho}\Big(\rho^{-2\alpha}G_\alpha\Big(\frac{r}{\rho}\Big)\Big)\\
&=(f(\rho)-f(0))^2\lim_{r\rightarrow0^+}\Big[-\frac{2\alpha}{\rho^{1+2\alpha}}G_\alpha\Big(\frac{r}{\rho}\Big)-\frac{r}{\rho^{2+2\alpha}}G'_\alpha\Big(\frac{r}{\rho}\Big)\Big]\\
&=-\frac{2\alpha G_\alpha(0)}{\rho^{1+2\alpha}}(f(\rho)-f(0))^2\\
&=-2\alpha B\Big(\frac12,\frac{n+1}{2}\Big)\frac{(f(\rho)-f(0))^2}{\rho^{1+2\alpha}}.
\end{split}
\end{eqnarray}
Thus,
\begin{eqnarray}\label{main-term-1}
\begin{split}
&\frac{c''_{n,\alpha}}{2}\int_0^\infty\Big(\frac{1}{r^{n+2\alpha}}
\frac{\partial}{\partial\rho}\Big(\rho^nG_\alpha\Big(\frac{\rho}{r}\Big)\Big)(f(\rho)-f(r))^2\Big)\Big|^{\infty}_{r\rightarrow0^+}d\rho\\
&=\alpha c''_{n,\alpha} B\Big(\frac12,\frac{n+1}{2}\Big)\int_0^\infty\frac{(f(\rho)-f(0))^2}{\rho^{1+2\alpha}}d\rho\\
&=\frac{\alpha c''_{n,\alpha}B(\frac12,\frac{n+1}{2})}{\omega_{n-1}}
\int_{\R^n}\frac{(f(0)-f(x))^2}{|x|^{n+2\alpha}}dx\\
&=\frac{\alpha2^{2\alpha-1}\Gamma(\frac{n}{2}+\alpha)}{\Gamma(1-\alpha)\Gamma(\frac{n}{2}+1)}
\int_{\R^n}\frac{(f(0)-f(x))^2}{|x|^{n+2\alpha}}dx.
\end{split}
\end{eqnarray}
We proceed to check that, for all $0<\rho,r<\infty$ and $\rho\neq r$,
\begin{eqnarray}\label{negative}
\frac{\partial}{\partial r}\Big(\frac{1}{r^{n+2\alpha}}\frac{\partial}{\partial\rho}\Big(\rho^nG_\alpha\Big(\frac{\rho}{r}\Big)\Big)\Big)
<0.
\end{eqnarray}
Indeed, direct calculations give that, for $0<\rho<r$,
\begin{eqnarray*}
\frac{\partial}{\partial\rho}
\Big(\rho^nG_\alpha\Big(\frac{\rho}{r}\Big)\Big)
=n\rho^{n-1}G_\alpha\Big(\frac{\rho}{r}\Big)
+\frac{\rho^n}{r}G'_\alpha\Big(\frac{\rho}{r}\Big)
\end{eqnarray*}
and
\begin{eqnarray}\label{negative-1}
\begin{split}
\frac{\partial}{\partial r}\Big(\frac{1}{r^{n+2\alpha}}\frac{\partial}{\partial\rho}\Big(\rho^nG_\alpha\Big(\frac{\rho}{r}\Big)\Big)\Big)
&=n\rho^{n-1}\frac{\partial}{\partial r}\Big(\frac{1}{r^{n+2\alpha}}G_\alpha\Big(\frac{\rho}{r}\Big)\Big)
+\rho^{n}\frac{\partial}{\partial r}\Big(\frac{1}{r^{n+1+2\alpha}}G'_\alpha\Big(\frac{\rho}{r}\Big)\Big)\\
&=-\frac{n(n+2\alpha)\rho^{n-1}}{r^{n+1+2\alpha}}G_\alpha\Big(\frac{\rho}{r}\Big)
-\frac{(2n+1+2\alpha)\rho^{n}}{r^{n+2+2\alpha}}G'_\alpha\Big(\frac{\rho}{r}\Big)\\
&\ \ \ \ \ \
-\frac{\rho^{n+1}}{r^{n+3+2\alpha}}G''_\alpha\Big(\frac{\rho}{r}\Big)<0,
\end{split}
\end{eqnarray}
where in the last inequality we have used Corollary \ref{Positivity-Convergence} showing that $G_\alpha$, $G'_\alpha$ and $G''_\alpha$ are all positive in $(0,1)$.

Furthermore, for $0<r<\rho$, by \eqref{transform-identity} and Corollary \ref{Positivity-Convergence}, we derive that
\begin{eqnarray}\label{negative-2}
\begin{split}
\frac{\partial}{\partial r}\Big(\frac{1}{r^{n+2\alpha}}\frac{\partial}{\partial\rho}\Big(\rho^{n}G_\alpha\Big(\frac{\rho}{r}\Big)\Big)\Big)
&=\frac{\partial}{\partial r}\Big(\frac{\partial}{\partial\rho}\Big(\rho^{-2\alpha}G\Big(\frac{r}{\rho}\Big)\Big)\Big)\\
&=\frac{\partial}{\partial r}\Big(-\frac{2\alpha}{\rho^{1+2\alpha}}G_\alpha\Big(\frac{r}{\rho}\Big)-\frac{r}{\rho^{2+2\alpha}}G'_\alpha\Big(\frac{r}{\rho}\Big)\Big)\\
&=-\frac{1+2\alpha}{\rho^{2+2\alpha}}G'_\alpha\Big(\frac{r}{\rho}\Big)
-\frac{r}{\rho^{3+2\alpha}}G''_\alpha\Big(\frac{r}{\rho}\Big)<0.
\end{split}
\end{eqnarray}
Combining \eqref{negative-1} and \eqref{negative-2} gives \eqref{negative}, which implies that
\begin{eqnarray*}
\int_0^\infty\int_0^\infty
\frac{\partial}{\partial r}\Big(\frac{1}{r^{n+2\alpha}}\frac{\partial}{\partial\rho}\Big(\rho^nG_\alpha\Big(\frac{\rho}{r}\Big)\Big)\Big)(f(\rho)-f(r))^2d\rho dr\leq0.
\end{eqnarray*}
This negative integral along with
\eqref{representation-1} and \eqref{main-term-1} completes the proof of Proposition \ref{n-nonlinear-inequality}. \hfill\hfill$\square$\vskip12pt
\begin{remark}
We remark that Li and Rodrigo \cite{[Li-Rodrigo20]} established the following weighted nonlinear inequality: for any $\delta\in(-2\alpha,2-2\alpha)$ and radial decreasing Schwartz function $f$,
\begin{eqnarray}\label{nonlinear-inequality}
\int_{\R^n}\frac{\Lambda^{-2+2\alpha}\nabla f(x)\cdot\nabla f(x)}{|x|^{n+\delta}}dx
\geq
C_{n,\alpha,\delta}\int_{\R^n}\frac{(f(0)-f(x))^2}{|x|^{n+2\alpha+\delta}}dx,
\end{eqnarray}
which can be applied to show the finite-time blow-up of smooth solutions to \eqref{M-CCF-T} for some class of radial decreasing initial data. The main technique in the proof of \eqref{nonlinear-inequality} from \cite{[Li-Rodrigo20]} is Hardy's inequality and a pointwise lower bound for the nonlinear term, i.e., for the radial monotone function $f$ and for any $x\in\mathbb{R}^n\backslash\{0\}$,
\begin{eqnarray*}
\Lambda^{-2+2\alpha}\nabla f(x)\cdot\nabla f(x)
\geq C_{n,\alpha}\frac{f'(|x|)}{|x|^{n-1+2\alpha}}
\int_0^{|x|} f'(\rho)\rho^nd\rho,
\end{eqnarray*}
which can be seen from the first equality in \eqref{expression-nonlinear-term} and the fact that $G_\alpha(\lambda)>G_\alpha(0)$ for $\lambda\in(0,1)$.
Proposition \ref{n-nonlinear-inequality} is the special case $\delta=0$ in \eqref{nonlinear-inequality}, but whose proof provided here does not depend on the monotone decaying property of $f$. Unfortunately, from the second equality in \eqref{condition}, we can see that the current method may not be applied to prove more general inequality \eqref{nonlinear-inequality} for  $\delta\in(0,2-2\alpha)$, which is needed for the implication of blow-up. Finally, one can refer to \cite{[Dong]} for another proof of \eqref{nonlinear-inequality} with the use of Meillin transform
for a general radial function $f$.
\end{remark}
Since the nonlinear inequality in Proposition \ref{n-nonlinear-inequality}, which involves a non-integrable weight $\frac{1}{|x|^n}$, is not directly useful for implying blow-ups of \eqref{M-CCF-T} for general radial initial data, we establish the next weighed inequality to show the finite-time blow-up. It should be remarked that the one-dimensional analogous inequality was proved by Li and Rodrigo in \cite{[Li-Rodrigo20]}.
\begin{proposition}\label{nonlinear-weighted-inequality-exponential}
Let $n\geq2$ and $0<\alpha<1$.
Let $f:\R^n\rightarrow\R$ be a radial Schwartz function. Then
\begin{eqnarray}\label{nonlinear-inequality-exponential-weight}
\int_{\R^n}\frac{\Lambda^{-2+2\alpha}\nabla f(x)\cdot\nabla f(x)}{|x|^{n}}e^{-|x|}dx
\geq C'_{n,\alpha}
\int_{\R^n}\frac{(f(0)-f(x))^2}{|x|^{n+2\alpha}}dx
-C''_{n,\alpha}\|f\|^2_{L^\infty},
\end{eqnarray}
where
\begin{eqnarray*}
C'_{n,\alpha}=\frac{\alpha2^{2\alpha-2}\Gamma(\frac{n}{2}+\alpha)}{\Gamma(1-\alpha)\Gamma(\frac{n}{2}+1)}
\end{eqnarray*}
and the constant $C''_{n,\alpha}$ depends only on $n$ and $\alpha$.
\end{proposition}
\textbf{Proof.} We only need to modify the proof of Proposition \ref{n-nonlinear-inequality}.
Similar to \eqref{representation-1}, by \eqref{expression-nonlinear-term}, \eqref{boundary-1}, \eqref{condition} and \eqref{negative},
we can derive that
\begin{eqnarray}\label{representation-2}
\begin{split}
&\int_{\R^n}\frac{\Lambda^{-2+2\alpha}\nabla f(x)\cdot\nabla f(x)}{|x|^{n}}e^{-|x|}dx\\
&=c''_{n,\alpha}\int_0^\infty\frac{f'(r)e^{-r}}{r^{n+2\alpha}}
\Big(-\int_0^\infty \frac{\partial}{\partial\rho}\Big(\rho^{n}G_\alpha\Big(\frac{\rho}{r}\Big)\Big)(f(\rho)-f(r))d\rho\Big)dr\\
&=\frac{c''_{n,\alpha}}{2}\int_0^\infty\int_0^\infty\frac{e^{-r}}{r^{n+2\alpha}}
 \frac{\partial}{\partial\rho}\Big(\rho^{n}G_\alpha\Big(\frac{\rho}{r}\Big)\Big)\frac{\partial}{\partial r}\Big((f(\rho)-f(r))^2\Big)drd\rho\\
&=\frac{c''_{n,\alpha}}{2}\int_0^\infty\Big(\frac{e^{-r}}{r^{n+2\alpha}}\frac{\partial}{\partial\rho}\Big(\rho^{n}G_\alpha\Big(\frac{\rho}{r}\Big)\Big)(f(\rho)-f(r))^2\Big)\Big|^{\infty}_{r\rightarrow0^+}d\rho\\
&\ \ \
-\frac{c''_{n,\alpha}}{2}\int_0^\infty\int_0^\infty
\frac{\partial}{\partial r}\Big(\frac{e^{-r}}{r^{n+2\alpha}}\frac{\partial}{\partial\rho}\Big(\rho^{n}G_\alpha\Big(\frac{\rho}{r}\Big)\Big)\Big)(f(\rho)-f(r))^2d\rho dr\\
&=\alpha c''_{n,\alpha} B\Big(\frac12,\frac{n+1}{2}\Big)\int_0^\infty\frac{(f(\rho)-f(0))^2}{\rho^{1+2\alpha}}d\rho\\
&\ \ \ \
-\frac{c''_{n,\alpha}}{2}\int_0^\infty\int_0^\infty
e^{-r}\frac{\partial}{\partial r}\Big(\frac{1}{r^{n+2\alpha}}\frac{\partial}{\partial\rho}\Big(\rho^{n}G_\alpha\Big(\frac{\rho}{r}\Big)\Big)\Big)(f(\rho)-f(r))^2d\rho dr\\
&\ \ \ \
+\frac{c''_{n,\alpha}}{2}\int_0^\infty\int_0^\infty
\frac{e^{-r}}{r^{n+2\alpha}}\frac{\partial}{\partial\rho}\Big(\rho^{n}G_\alpha\Big(\frac{\rho}{r}\Big)\Big)\Big)(f(\rho)-f(r))^2d\rho dr\\
&\geq
\frac{\alpha2^{2\alpha-1}\Gamma(\frac{n}{2}+\alpha)}{\Gamma(1-\alpha)\Gamma(\frac{n}{2}+1)}
\int_{\R^n}\frac{(f(0)-f(x))^2}{|x|^{n+2\alpha}}dx\\
& \ \  \
+\frac{c''_{n,\alpha}}{2}\underbrace{\int_0^\infty\int_0^\infty
\frac{e^{-r}}{r^{n+2\alpha}}\frac{\partial}{\partial\rho}\Big(\rho^{n}G_\alpha\Big(\frac{\rho}{r}\Big)\Big)(f(\rho)-f(r))^2d\rho dr}_{\mathcal{I}},
\end{split}
\end{eqnarray}
where the constant $c''_{n,\alpha}$ is defined as \eqref{constant-11}.

We proceed to estimate the integral $\mathcal{I}$.
For $0<\rho<r$, straightforward computation yields that
\begin{eqnarray*}
\frac{1}{r^{n+2\alpha}}\frac{\partial}{\partial\rho}\Big(\rho^{n}G_\alpha\Big(\frac{\rho}{r}\Big)\Big)\Big)
=\frac{\rho^{n-1}}{r^{n+2\alpha}}\Big[nG_\alpha\Big(\frac{\rho}{r}\Big)+\frac{\rho}{r}G'_\alpha\Big(\frac{\rho}{r}\Big)\Big].
\end{eqnarray*}
For $0<r<\rho$, by \eqref{transform-identity}, we have that
\begin{align*}
\frac{1}{r^{n+2\alpha}}\frac{\partial}{\partial\rho}\Big(\rho^{n}G_\alpha\Big(\frac{\rho}{r}\Big)\Big)\Big)
&=\frac{\partial}{\partial\rho}\Big(\rho^{-2\alpha}G_\alpha\Big(\frac{r}{\rho}\Big)\Big)\Big)\\
&=-\frac{2\alpha}{\rho^{1+2\alpha}}G_\alpha\Big(\frac{r}{\rho}\Big)
-\frac{r}{\rho^{2+2\alpha}}G'_\alpha\Big(\frac{r}{\rho}\Big).
\end{align*}
Thus, we rewrite the integral $\mathcal{I}$ as
\begin{eqnarray}\label{potential-negative-term}
\begin{split}
\mathcal{I}&=\iint\limits_{0<\rho<r}\frac{e^{-r}\rho^{n-1}(f(\rho)-f(r))^2}{r^{n+2\alpha}}\Big[nG_\alpha\Big(\frac{\rho}{r}\Big)+\frac{\rho}{r}G'_\alpha\Big(\frac{\rho}{r}\Big)\Big]d\rho dr\\
&\ \ \ \ \
-\iint\limits_{0<r<\rho}\frac{e^{-r}(f(\rho)-f(r))^2}{\rho^{1+2\alpha}}\Big[2\alpha G_\alpha\Big(\frac{r}{\rho}\Big)
+\frac{r}{\rho}G'_\alpha\Big(\frac{r}{\rho}\Big)\Big]d\rho dr.
\end{split}
\end{eqnarray}
Next we will divide into three cases $0<\alpha<\frac12$, $\frac{1}{2}<\alpha<1$ and $\alpha=\frac12$ to derive a lower bound for the integral $\mathcal{I}$.

\textbf{Case 1}. $0<\alpha<\frac12$. In this case, by Corollary \ref{Positivity-Convergence}, \eqref{G-definition} and \eqref{1-order} leads to
\begin{align*}
\mathcal{I}
&\geq-\iint\limits_{0<r<\rho}\frac{e^{-r}(f(\rho)-f(r))^2}{\rho^{1+2\alpha}}\Big[2\alpha G_\alpha\Big(\frac{r}{\rho}\Big)
+\frac{r}{\rho}G'_\alpha\Big(\frac{r}{\rho}\Big)\Big]d\rho dr\\
&=-\iint\limits_{0<r<\rho}\frac{e^{-r}(f(\rho)-f(r))^2}{\rho^{1+2\alpha}}
\Big[\sum^\infty_{k=0}2(k+\alpha)a_{2k}(\alpha)\Big(\frac{r}{\rho}\Big)^{2k}\Big]d\rho dr\\
&=-\sum^\infty_{k=0}2(k+\alpha)a_{2k}(\alpha)\iint\limits_{0<r<\rho}\frac{e^{-r}(f(\rho)-f(r))^2}{\rho^{1+2\alpha}}
\Big(\frac{r}{\rho}\Big)^{2k}d\rho dr\\
&\geq-4\|f\|^2_{L^\infty}\sum^\infty_{k=0}2(k+\alpha)a_{2k}(\alpha)\int^\infty_0e^{-r}r^{2k}\Big(\int^\infty_r\frac{d\rho}{\rho^{1+2\alpha+2k}}\Big)dr\\
&=-4\Gamma(1-2\alpha)\cdot\sum^\infty_{k=0}a_{2k}(\alpha)\cdot\|f\|^2_{L^\infty},
\end{align*}
which along with \eqref{representation-2} implies the desired inequality \eqref{nonlinear-inequality-exponential-weight}.

\textbf{Case 2}. $\frac{1}{2}<\alpha<1$. In this case, by \eqref{potential-negative-term}, Corollary \ref{Positivity-Convergence} and a change of variable,
we derive that
\begin{eqnarray}\label{3.13}
\begin{split}
\mathcal{I}&\geq\iint\limits_{0<\rho<r}\frac{e^{-r}\rho^{n-1}}{r^{n+2\alpha}}(f(\rho)-f(r))^2\frac{\rho}{r}G'_\alpha\Big(\frac{\rho}{r}\Big)d\rho dr\\
&\ \ \ \ \
-\iint\limits_{0<r<\rho}\frac{e^{-r}(f(\rho)-f(r))^2}{\rho^{1+2\alpha}}\Big[2\alpha G_\alpha\Big(\frac{r}{\rho}\Big)
+\frac{r}{\rho}G'_\alpha\Big(\frac{r}{\rho}\Big)\Big]d\rho dr\\
&=-\underbrace{\iint\limits_{0<\rho<r}\Big[e^{-\rho}-\frac{e^{-r}\rho^{n-1}}{r^{n-1}}\Big]
\frac{(f(\rho)-f(r))^2}{r^{1+2\alpha}}\frac{\rho}{r}G'_\alpha\Big(\frac{\rho}{r}\Big)d\rho dr}_{\mathcal{J}}\\
&\ \ \ \ \
-2\alpha\underbrace{\iint\limits_{0<r<\rho}\frac{e^{-r}(f(\rho)-f(r))^2}{\rho^{1+2\alpha}} G_\alpha\Big(\frac{r}{\rho}\Big)
d\rho dr}_{\mathcal{K}}.
\end{split}
\end{eqnarray}
By \eqref{G-definition}, $\frac12<\alpha<1$ and Corollary \ref{Positivity-Convergence},  the positive term $\mathcal{K}$ can be estimated above as can be estimated above as
\begin{eqnarray}\label{3.14}
\begin{split}
\mathcal{K}
&=\sum^\infty_{k=0}a_{2k}(\alpha)\iint\limits_{0<r<\rho}\frac{e^{-r}(f(\rho)-f(r))^2}{\rho^{1+2\alpha}} \Big(\frac{r}{\rho}\Big)^{2k}d\rho dr\\
&\leq\sum^\infty_{k=0}2a_{2k}(\alpha)\Big[\iint\limits_{0<r<\rho}\frac{(f(\rho)-f(0))^2}{\rho^{1+2\alpha+2k}} r^{2k}d\rho dr+\iint\limits_{0<r<\rho}\frac{(f(r)-f(0))^2}{\rho^{1+2\alpha+2k}}r^{2k} d\rho dr\Big]\\
 &=\sum^\infty_{k=0}a_{2k}(\alpha)\Big(\frac{2}{2k+1}+\frac{1}{k+\alpha}\Big)
 \int^\infty_0\frac{(f(\rho)-f(0))^2}{\rho^{2\alpha}}d\rho\\
 &\leq\sum^\infty_{k=0}\frac{4a_{2k}(\alpha)}{2k+1}\cdot
 \int^\infty_0\frac{(f(\rho)-f(0))^2}{\rho^{2\alpha}}d\rho.
\end{split}
\end{eqnarray}
By \eqref{1-order}, we further express the term $\mathcal{J}$ as
\begin{eqnarray}\label{J-term}
\begin{split}
\mathcal{J}
&=\iint\limits_{0<\rho<r}\Big[e^{-\rho}-e^{-r}\Big(\frac{\rho}{r}\Big)^{n-1}\Big]
\frac{(f(\rho)-f(r))^2}{r^{1+2\alpha}}\sum^\infty_{k=1}2ka_{2k}(\alpha)\Big(\frac{\rho}{r}\Big)^{2k}d\rho dr\\
&=\sum^\infty_{k=1}2ka_{2k}(\alpha)\underbrace{\iint\limits_{0<\rho<r}\Big[e^{-\rho}-e^{-r}\Big(\frac{\rho}{r}\Big)^{n-1}\Big]
\frac{(f(\rho)-f(r))^2}{r^{1+2\alpha}}\Big(\frac{\rho}{r}\Big)^{2k}d\rho dr}_{\mathcal{J}_1}.
\end{split}
\end{eqnarray}
By the mean-value theorem, we derive that, for $0<\rho<r$,
\begin{eqnarray}\label{mean-value}
\begin{split}
e^{-\rho}-e^{-r}\Big(\frac{\rho}{r}\Big)^{n-1}
&=\rho^{n-1}\Big(\frac{e^{-\rho}}{\rho^{n-1}}-\frac{e^{-r}}{r^{n-1}}\Big)\\
&\leq e^{-\rho}\rho^{n-1}
\Big(\rho^{1-n}+\frac{n-1}{\rho^{n}}\Big)(r-\rho)\\
&=e^{-\rho}(r-\rho)
+(n-1)e^{-\rho}\Big(\frac{r}{\rho}-1\Big).
\end{split}
\end{eqnarray}
Then $\mathcal{J}_1$ can be bounded as
\begin{eqnarray}\label{J1}
\begin{split}
\mathcal{J}_1
&\leq
\underbrace{\iint\limits_{0<\rho<r}
e^{-\rho}(r-\rho)\frac{(f(\rho)-f(r))^2}{r^{1+2\alpha}}\Big(\frac{\rho}{r}\Big)^{2k}d\rho dr}_{\mathcal{J}_{11}}\\
&\ \ \
+
(n-1)\underbrace{\iint\limits_{0<\rho<r}
e^{-\rho}\Big(\frac{r}{\rho}-1\Big)\frac{(f(\rho)-f(r))^2}{r^{1+2\alpha}}\Big(\frac{\rho}{r}\Big)^{2k}d\rho dr}_{\mathcal{J}_{12}}.
\end{split}
\end{eqnarray}
Since $\frac{1}{2}<\alpha<1$, we derive that, for any $k\geq1$,
\begin{eqnarray}\label{J11}
\begin{split}
\mathcal{J}_{11}
&\leq
4\|f\|^2_{L^\infty}\iint\limits_{0<\rho<r}
e^{-\rho}\frac{r-\rho}{r^{1+2\alpha}}\Big(\frac{\rho}{r}\Big)^{2k}d\rho
dr\\
&=4\|f\|^2_{L^\infty}\int^\infty_0e^{-\rho}\rho^{2k}\Big(\int^\infty_\rho\frac{1-\frac{\rho}{r}}{r^{2k+2\alpha}}dr\Big)
d\rho\\
&=4\|f\|^2_{L^\infty}\Big(\frac{1}{2k+2\alpha-1}
-\frac{1}{2k+2\alpha}\Big)\int^\infty_0e^{-\rho}\rho^{1-2\alpha}d\rho\\
&=\frac{4\Gamma(2-2\alpha)\|f\|^2_{L^\infty}}{(2k+2\alpha-1)(2k+2\alpha)}
<\frac{\Gamma(2-2\alpha)}{k^2}\|f\|^2_{L^\infty}.
\end{split}
\end{eqnarray}
Note that $\mathcal{J}_{12}$ can not be estimated in the way same as $\mathcal{J}_{11}$.
Utilizing the elementary inequality
\begin{eqnarray*}
(a-b)^2\leq2(a-c)^2+2(b-c)^2~{\rm  for}~ a,b,c\in\mathbb{R},
\end{eqnarray*}
we estimate $\mathcal{J}_{12}$ in the following way, for any $k\geq1$,
\begin{align*}
\mathcal{J}_{12}
&\leq2\iint\limits_{0<\rho<r}
(\frac{r}{\rho}-1)\frac{(f(\rho)-f(0))^2}{r^{1+2\alpha}}(\frac{\rho}{r})^{2k}d\rho dr
+
2\iint\limits_{0<\rho<r}
(\frac{r}{\rho}-1)\frac{(f(0)-f(r))^2}{r^{1+2\alpha}}(\frac{\rho}{r})^{2k}d\rho dr\\
&=2\int\limits_0^\infty\rho^{2k}(f(\rho)-f(0))^2(\int\limits^\infty_\rho\frac{\frac{r}{\rho}-1}{r^{2k+2\alpha+1}}dr)d\rho
+2\int\limits_0^\infty\frac{(f(r)-f(0))^2}{r^{2k+2\alpha+1}}(\int\limits^r_{0}\rho^{2k}(\frac{r}{\rho}-1)d\rho)dr\\
&=\frac{2}{(2k+2\alpha-1)(2k+2\alpha)}\int\limits^\infty_0\frac{(f(\rho)-f(0))^2}{\rho^{2\alpha}}d\rho
+\frac{1}{k(2k+1)}\int\limits_0^\infty\frac{(f(r)-f(0))^2}{r^{2\alpha}}dr\\
&\leq\frac{1}{k^2}\int^\infty_0\frac{(f(\rho)-f(0))^2}{\rho^{2\alpha}}d\rho,
\end{align*}
which follows from \eqref{J11}, \eqref{J1} and \eqref{J-term} that
\begin{eqnarray*}
\mathcal{J}
\leq2(n-1)\sum^\infty_{k=1}\frac{a_{2k}(\alpha)}{k}\int^\infty_0\frac{(f(\rho)-f(0))^2}{\rho^{2\alpha}}d\rho
+2\Gamma(2-2\alpha)\sum^\infty_{k=1}\frac{a_{2k}(\alpha)}{k}\cdot\|f\|^2_{L^\infty}.
\end{eqnarray*}
This estimate along with \eqref{3.13} and \eqref{3.14} yields that
\begin{align*}
\mathcal{I}
&\geq
-\Big[2(n-1)\sum^\infty_{k=1}\frac{a_{2k}(\alpha)}{k}
+8\alpha\sum^\infty_{k=0}\frac{a_{2k}(\alpha)}{2k+1}\Big]
\int^\infty_0\frac{(f(\rho)-f(0))^2}{\rho^{2\alpha}}d\rho\\
&\ \ \ \ \ \ \
-2\Gamma(2-2\alpha)\sum^\infty_{k=1}\frac{a_{2k}(\alpha)}{k}\cdot\|f\|^2_{L^\infty},
\end{align*}
Then, by utilizing \eqref{representation-2} and the first inequality in Lemma \ref{Young-inequality} with
\begin{eqnarray*}
\epsilon=\frac{\alpha B(\frac12,\frac{n+1}{2})}
{2(n-1)\displaystyle\sum^\infty_{k=1}\frac{a_{2k}(\alpha)}{k}
+8\alpha\displaystyle\sum^\infty_{k=0}\frac{a_{2k}(\alpha)}{2k+1}},
\end{eqnarray*}
we can obtain the desired inequality \eqref{nonlinear-inequality-exponential-weight}.

\textbf{Case 3}. $\alpha=\frac12$. In this case, for the conciseness, we use some abbreviations of the notation as follows,
\begin{eqnarray*}
G_{\frac12}(\lambda):=G(\lambda),~a_{2k}\Big(\frac12\Big):=a_{2k}
\end{eqnarray*}
for $\lambda\in(0,1)$ and $k=0,1,2\cdot\cdot\cdot$.
Similar to \eqref{3.13}, from \eqref{potential-negative-term}, we have that
\begin{eqnarray}\label{I-term-case-3}
\begin{split}
\mathcal{I}
&\geq
-\underbrace{\iint\limits_{0<\rho<r}\Big[e^{-\rho}-\frac{e^{-r}\rho^{n-1}}{r^{n-1}}\Big]
\frac{(f(\rho)-f(r))^2}{r^{2}}\frac{\rho}{r}G'\Big(\frac{\rho}{r}\Big)d\rho dr}_{\mathcal{\widetilde{J}}}\\
&\ \ \ \ \
-\underbrace{\iint\limits_{0<r<\rho}\frac{e^{-r}(f(\rho)-f(r))^2}{\rho^{2}} G\Big(\frac{r}{\rho}\Big)
d\rho dr}_{\mathcal{\widetilde{K}}}.
\end{split}
\end{eqnarray}
By \eqref{G-definition}, we split the term $\mathcal{\widetilde{K}}$ into two parts as follows
\begin{align*}
\mathcal{\widetilde{K}}
&=\sum^\infty_{k=0}a_{2k}\iint\limits_{0<r<\rho}\frac{e^{-r}(f(\rho)-f(r))^2}{\rho^{2}} \Big(\frac{r}{\rho}\Big)^{2k}
d\rho dr\\
&=\sum^\infty_{k=0}a_{2k}
\underbrace{\iint\limits_{0<r\leq\frac{\rho}{2}}\frac{(f(\rho)-f(r))^2}{e^{r}\rho^{2}} (\frac{r}{\rho})^{2k}d\rho dr}_{\mathcal{\widetilde{K}}_1}
+\sum^\infty_{k=0}a_{2k}
\underbrace{\iint\limits_{\frac{\rho}{2}<r<\rho}\frac{(f(\rho)-f(r))^2}{e^r\rho^{2}} (\frac{r}{\rho})^{2k}d\rho dr}_{\mathcal{\widetilde{K}}_2}.
\end{align*}
The term $\mathcal{\widetilde{K}}_1$ can be estimated as
\begin{align*}
\mathcal{\widetilde{K}}_1
&\leq\Big(\frac12\Big)^{2k}\iint\limits_{0<r<\frac{\rho}{2}}\frac{e^{-r}(f(\rho)-f(r))^2}{\rho^{2}}
d\rho dr\\
&\leq
\frac{1}{2^{2k-1}}\iint\limits_{0<r<\frac{\rho}{2}}\frac{(f(\rho)-f(0))^2}{\rho^{2}e^r}
d\rho dr+\frac{1}{2^{2k-1}}\iint\limits_{0<r<\frac{\rho}{2}}\frac{(f(r)-f(0))^2}{\rho^2e^r}
d\rho dr\\
&=\frac{1}{2^{2k-1}}\int_0^\infty\frac{(f(\rho)-f(0))^2}{\rho^{2}}(1-e^{-\frac{\rho}{2}})d\rho
+\frac{1}{2^{2k}}\int_0^\infty \frac{(f(r)-f(0))^2}{re^{r}}dr.
\end{align*}
For the term $\mathcal{\widetilde{K}}_2$, we derive that
\begin{align*}
\mathcal{\widetilde{K}}_2
&\leq2\iint\limits_{\frac{\rho}{2}<r<\rho}\frac{e^{-r}(f(\rho)-f(0))^2}{\rho^{2}} \Big(\frac{r}{\rho}\Big)^{2k}
d\rho dr+2\iint\limits_{\frac{\rho}{2}<r<\rho}\frac{e^{-r}(f(r)-f(0))^2}{\rho^{2}} \Big(\frac{r}{\rho}\Big)^{2k}
d\rho dr\\
&\leq2\int_0^\infty\frac{(f(\rho)-f(0))^2}{\rho^{2+2k}}\cdot \frac{e^{-\frac{\rho}{2}}\rho^{2k+1}}{2k+1}d\rho
+2\int_0^\infty \frac{r^{2k}(f(r)-f(0))^2}{e^{r}}\cdot\frac{r^{-2k-1}}{2k+1}dr\\
&=\frac{2}{2k+1}\int_0^\infty\frac{(f(\rho)-f(0))^2}{\rho e^{\frac\rho2}}d\rho
+\frac{2}{2k+1}\int_0^\infty \frac{(f(r)-f(0))^2}{re^{r}}dr\\
&\leq\frac{4}{2k+1}\int_0^\infty\frac{(f(\rho)-f(0))^2}{\rho e^{\frac\rho2}}d\rho.
\end{align*}
Thus, we obtain that
\begin{eqnarray}\label{K-term-Case-3}
\begin{split}
\mathcal{\widetilde{K}}
&\leq\sum^\infty_{k=0}\frac{a_{2k}}{2^{2k-1}}\int\limits_0^\infty\frac{(f(\rho)-f(0))^2}{\rho^{2}}(1-e^{-\frac{\rho}{2}})d\rho
+\sum^\infty_{k=0}(\frac{a_{2k}}{2^{2k}}+\frac{4a_{2k}}{2k+1})\int\limits_0^\infty\frac{(f(\rho)-f(0))^2}{\rho e^{\frac\rho2}}d\rho\\
&\leq\sum^\infty_{k=0}\frac{2a_{2k}}{2k+1}\int\limits_0^\infty\frac{(f(\rho)-f(0))^2}{\rho^{2}}(1-e^{-\frac{\rho}{2}})d\rho
+\sum^\infty_{k=0}\frac{5a_{2k}}{2k+1}\int\limits_0^\infty\frac{(f(\rho)-f(0))^2}{\rho e^{\frac\rho2}}d\rho.
\end{split}
\end{eqnarray}
From \eqref{1-order} and \eqref{mean-value}, we can bound $\mathcal{\widetilde{J}}$ as
\begin{align*}
\mathcal{\widetilde{J}}
&=\sum^\infty_{k=1}2ka_{2k}
\iint\limits_{0<\rho<r}\Big[e^{-\rho}-e^{-r}\Big(\frac{\rho}{r}\Big)^{n-1}\Big]
\frac{(f(\rho)-f(r))^2}{r^2}\Big(\frac{\rho}{r}\Big)^{2k}d\rho dr\\
&\leq\sum^\infty_{k=1}2ka_{2k}
\underbrace{\iint\limits_{0<\rho<r}
e^{-\rho}(r-\rho)\frac{(f(\rho)-f(r))^2}{r^2}\Big(\frac{\rho}{r}\Big)^{2k}d\rho dr}_{\mathcal{\widetilde{J}}_{1}}\\
&\
+(n-1)\sum^\infty_{k=1}2ka_{2k}
\underbrace{\iint\limits_{0<\rho<r}
e^{-\rho}\Big(\frac{r}{\rho}-1\Big)\frac{(f(\rho)-f(r))^2}{r^2}\Big(\frac{\rho}{r}\Big)^{2k}d\rho dr}_{\mathcal{\widetilde{J}}_{2}}.
\end{align*}
Similar to \eqref{J11}, $\mathcal{\widetilde{J}}_1$ can be estimated as, for any $k\geq1$,
\begin{eqnarray}\label{J1-Case-3}
\begin{split}
\mathcal{\widetilde{J}}_1
&\leq
4\|f\|^2_{L^\infty}\iint\limits_{0<\rho<r}
e^{-\rho}\frac{r-\rho}{r^2}\Big(\frac{\rho}{r}\Big)^{2k}d\rho
dr\\
&=4\|f\|^2_{L^\infty}\int^\infty_0e^{-\rho}\rho^{2k}\Big(\int^\infty_\rho\frac{1-\frac{\rho}{r}}{r^{2k+1}}dr\Big)
d\rho\\
&=4\|f\|^2_{L^\infty}\Big(\frac{1}{2k}
-\frac{1}{2k+1}\Big)\int^\infty_0e^{-\rho}d\rho\\
&=\frac{2\|f\|^2_{L^\infty}}{k(2k+1)}<\frac{1}{k^2}\|f\|^2_{L^\infty}.
\end{split}
\end{eqnarray}
We split the term $\mathcal{\widetilde{J}}_2$ into two parts in the following way
\begin{eqnarray*}
\mathcal{\widetilde{J}}_2
=\underbrace{\iint\limits_{0<\rho\leq\frac{r}{2}}
(\frac{r}{\rho}-1)\frac{(f(\rho)-f(r))^2}{r^2e^\rho}(\frac{\rho}{r})^{2k}d\rho dr}_{\mathcal{\widetilde{J}}_{21}}
+\underbrace{\iint\limits_{\frac{r}{2}<\rho<r}
(\frac{r}{\rho}-1)\frac{(f(\rho)-f(r))^2}{r^2e^\rho}(\frac{\rho}{r})^{2k}d\rho dr}_{\mathcal{\widetilde{J}}_{22}}.
\end{eqnarray*}
We estimate the positive term $\mathcal{\widetilde{J}}_{21}$ as, for any $k\geq1$,
\begin{align*}
\mathcal{\widetilde{J}}_{21}
&\leq\Big(\frac{1}{2}\Big)^{2k-1}
\iint\limits_{0<\rho\leq\frac{r}{2}}
e^{-\rho}\Big(\frac{r}{\rho}-1\Big)\frac{(f(\rho)-f(r))^2}{r^{2}}\cdot\frac{\rho}{r}d\rho dr\\
&=\Big(\frac{1}{2}\Big)^{2k-1}
\iint\limits_{0<\rho\leq\frac{r}{2}}
\Big(1-\frac{\rho}{r}\Big)\frac{(f(\rho)-f(r))^2}{r^2e^\rho}d\rho dr
\leq\frac{1}{k^2}
\iint\limits_{0<\rho\leq\frac{r}{2}}
\frac{(f(\rho)-f(r))^2}{r^2e^\rho}d\rho dr\nonumber\\
&\leq\frac{2}{k^2}
\iint\limits_{0<\rho\leq\frac{r}{2}}
\frac{(f(\rho)-f(0))^2}{r^2e^\rho}d\rho dr
+\frac{2}{k^2}
\iint\limits_{0<\rho\leq\frac{r}{2}}
\frac{(f(0)-f(r))^2}{r^2e^\rho}d\rho dr\nonumber\\
&=\frac{1}{k^2}\int^\infty_0\frac{(f(\rho)-f(0))^2}{\rho e^\rho}d\rho
+\frac{2}{k^2}\int^\infty_0\frac{(f(0)-f(r))^2}{r^{2}}(1-e^{-\frac{r}{2}}) dr.
\end{align*}
We further estimate the term $\mathcal{\widetilde{J}}_{22}$ in the following way, We further estimate the term $\mathcal{\widetilde{J}}_{22}$ in the following way,
\begin{align*}
\mathcal{\widetilde{J}}_{22}
&\leq2\iint\limits_{\frac{r}{2}<\rho<r}
(\frac{r}{\rho}-1)\frac{(f(\rho)-f(0))^2}{r^2e^{\frac{\rho}{2}}}(\frac{\rho}{r})^{2k}d\rho dr
+2\iint\limits_{\frac{r}{2}<\rho<r}
(\frac{r}{\rho}-1)\frac{(f(0)-f(r))^2}{r^2e^{\frac{r}{2}}}(\frac{\rho}{r})^{2k}d\rho dr\\
&=2\int\limits_0^\infty \frac{\rho^{2k}(f(\rho)-f(0))^2}{e^{\frac{\rho}{2}}}\Big[\int\limits^{2\rho}_\rho\frac{\frac{r}{\rho}-1}{r^{2k+2}}dr\Big]d\rho
+2\int\limits_0^\infty \frac{(f(r)-f(0))^2}{r^{2k+2}e^{\frac{r}{2}}}\Big[\int\limits^r_{\frac{r}{2}}\rho^{2k}(\frac{r}{\rho}-1)d\rho\Big]dr\\
&=4\Big[\frac{1-2^{-2k}}{2k}-\frac{1-2^{-2k-1}}{2k+1}\Big]\int^\infty_0\frac{(f(\rho)-f(0))^2}{\rho e^{\frac{\rho}{2}}}d\rho\\
&\leq4\Big[\frac{1}{2k}-\frac{1}{2k+1}\Big]\int^\infty_0\frac{(f(\rho)-f(0))^2}{\rho e^{\frac{\rho}{2}}}d\rho
\leq\frac{1}{k^2}\int^\infty_0\frac{(f(\rho)-f(0))^2}{\rho e^{\frac{\rho}{2}}}d\rho.
\end{align*}
Thus, we obtain that, for any $k\geq1$,
\begin{eqnarray*}
\mathcal{\widetilde{J}}_2
\leq
\frac{2}{k^2}\int^\infty_0\frac{(f(\rho)-f(0))^2}{\rho e^{\frac{\rho}{2}}}d\rho
+\frac{2}{k^2}\int^\infty_0\frac{(f(0)-f(r))^2}{r^{2}}(1-e^{-\frac{r}{2}}) dr.
\end{eqnarray*}
which along with \eqref{J1-Case-3} implies that
\begin{align*}
\mathcal{\widetilde{J}}
&\leq
4(n-1)\sum^\infty_{k=1}\frac{a_{2k}}{k}\cdot\int^\infty_0\frac{(f(\rho)-f(0))^2}{\rho e^{\frac{\rho}{2}}}d\rho
+
2\sum^\infty_{k=1}\frac{a_{2k}}{k}\cdot\|f\|^2_{L^\infty}\\
&\ \ \
+4(n-1)\sum^\infty_{k=1}\frac{a_{2k}}{k}\cdot\int^\infty_0\frac{(f(0)-f(r))^2}{r^{2}}(1-e^{-\frac{r}{2}}) dr.
\end{align*}
It follows from \eqref{K-term-Case-3} and \eqref{I-term-case-3} that
\begin{align*}
\mathcal{I}
&\geq
-\Big[4(n-1)\sum^\infty_{k=1}\frac{a_{2k}}{k}+\sum^\infty_{k=0}\frac{5a_{2k}}{2k+1}\Big]
\int^\infty_0\frac{(f(\rho)-f(0))^2}{\rho e^{\frac{\rho}{2}}}d\rho-2\sum^\infty_{k=1}\frac{a_{2k}}{k}\cdot\|f\|^2_{L^\infty}\\
&\ \ \
-\Big[4(n-1)\sum^\infty_{k=1}\frac{a_{2k}}{k}+\sum^\infty_{k=0}\frac{2a_{2k}}{2k+1}\Big]\int_0^\infty\frac{(f(\rho)-f(0))^2}{\rho^{2}}(1-e^{-\frac{\rho}{2}})d\rho.
\end{align*}
Taking
\begin{eqnarray*}
\epsilon=\frac{\alpha B(\frac12,\frac{n+1}{2})}{2\Big[4(n-1)\displaystyle\sum^\infty_{k=1}\frac{a_{2k}}{k}+\displaystyle\sum^\infty_{k=0}\frac{5a_{2k}}{2k+1}\Big]}
\end{eqnarray*}
in the second inequality of Lemma \ref{Young-inequality} and  another
\begin{eqnarray*}
\epsilon=\frac{\alpha B(\frac12,\frac{n+1}{2})}{2\Big[4(n-1)\displaystyle\sum^\infty_{k=1}\frac{a_{2k}}{k}+\displaystyle\sum^\infty_{k=0}\frac{2a_{2k}}{2k+1}\Big]}
\end{eqnarray*}
in the third inequality of Lemma \ref{Young-inequality}, it follows from \eqref{representation-2} that the desired inequality \eqref{nonlinear-inequality-exponential-weight}.
We then complete the proof of Proposition \ref{nonlinear-weighted-inequality-exponential}. \hfill\hfill$\square$\vskip12pt
\begin{remark}
We note that the case $\alpha=\frac12$ of \eqref{nonlinear-inequality-exponential-weight} in Proposition \ref{nonlinear-weighted-inequality-exponential} provides an inequality for the Riesz transform $\mathcal{R}:=-\Lambda^{-1}\nabla$:
for any radial Schwartz function $f:\R^n\rightarrow\R$, we have
\begin{eqnarray*}
-\int_{\R^n}\frac{\mathcal{R} f(x)\cdot\nabla f(x)}{|x|^{n}}e^{-|x|}dx
\geq \frac{1}{2nB(\frac12,\frac{n}{2})}
\int_{\R^n}\frac{(f(0)-f(x))^2}{|x|^{n+1}}dx
-C_{n}\|f\|^2_{L^\infty},
\end{eqnarray*}
where the constant $C_{n}$ depends only on $n$.
This inequality is a multi-dimensional generalization of the similar one for Hilbert transform $H:=-\Lambda^{-1}\partial_x$.
\end{remark}
\section{Proof of Theorem \ref{the-1}}
With Proposition \ref{nonlinear-weighted-inequality-exponential} in hand, we are now ready to prove Theorem \ref{the-1}.

\textbf{Proof of Theorem \ref{the-1}}. We will argue by contradiction.
Assume that the solution $\theta$ to \eqref{M-CCF-T} starting from the initial data $\theta_0\in \mathcal{S}(\R^n)$, the Schwartz class, satisfying \eqref{Class-initial-data} exists for all time. For our purpose, we introduce a quantity $J(t)$ given by
\begin{eqnarray}\label{weighted-integral-solution}
J(t)\triangleq\int_{\R^n}\frac{\theta(0,t)-\theta(x,t)}{|x|^{n}}e^{-|x|}dx.
\end{eqnarray}
By H\"{o}lder's inequality and Lemma \ref{local-well-posedness}, we have that
\begin{align*}
|J(t)|
&\leq \|\nabla\theta(t)\|_{L^\infty}\int_{|x|\leq1}\frac{dx}{|x|^{n-1}}
+2\|\theta_0\|_{L^\infty}\int_{|x|>1}\frac{e^{-|x|}}{|x|^{n}}dx\\
&=\omega_{n-1}\|\nabla\theta\|_{L^\infty}+2\omega_{n-1}\Big(\int^\infty_1\frac{r}{e^r}dr\Big)\|\theta_0\|_{L^\infty}<+\infty,
\end{align*}
which shows that $J(t)$ is finite for any $t>0$. Next we prove that $J(t)$ will blow up at some finite time $T_0>0$ and then obtain a contradiction.

By Lemma \ref{radial-property-preserved}, we know that the velocity at the origin is $0$, that is,
\begin{align*}
u(0,t)&=c_{n,\alpha}P.V.\int_{\R^n}\frac{\nabla\theta(x,t)}{|x|^{n-2+2\alpha}}dx\\
&=c_{n,\alpha}P.V.\int_{\R^n}\frac{x\theta'(|x|,t)}{|x|^{n-1+2\alpha}}dx=0,
\end{align*}
which along with \eqref{M-CCF-T}, Proposition \ref{nonlinear-weighted-inequality-exponential} and the maximum principle  $\|\theta(t)\|_{L^\infty}\leq\|\theta_0\|_{L^\infty}$ implies that
\begin{align*}
\frac{d}{dt}J(t)
&=\int_{\R^n} \frac{\Lambda^{-2+2\alpha}\nabla\theta(x,t)\cdot\nabla\theta(x,t)}{|x|^{n}}e^{-|x|}dx\\
&\geq C'_{n,\alpha}\int_{\R^n}\frac{(\theta(0,t)-\theta(x,t))^2}{|x|^{n+2\alpha}}dx
-C''_{n,\alpha}\|\theta_0\|^2_{L^\infty}.
\end{align*}
By \eqref{weighted-integral-solution} and H\"{o}lder's inequality, we note that
\begin{align*}
|J(t)|
&\leq\Big(\int_{\R^n}\frac{|\theta(0,t)-\theta(x,t)|^2}{|x|^{n+2\alpha}}dx\Big)^{\frac12}
\Big(\int_{\R^n}\frac{e^{-2|x|}}{|x|^{n-2\alpha}}dx\Big)^{\frac12}\\
&=\frac{\omega_{n-1}\Gamma(2\alpha)}{2^{2\alpha}}\Big(\int_{\R^n}\frac{|\theta(0,t)-\theta(x,t)|^2}{|x|^{n+2\alpha}}dx\Big)^{\frac12},
\end{align*}
which gives that
\begin{eqnarray*}
\int_{\R^n}\frac{|\theta(0,t)-\theta(x,t)|^2}{|x|^{n+2\alpha}}dx
\geq
\frac{2^{4\alpha}[J(t)]^2}{\omega^2_{n-1}\Gamma^2(2\alpha)}.
\end{eqnarray*}
Therefore, we derive that
\begin{eqnarray}\label{ordinary-differential-inequality}
\begin{split}
\frac{d}{dt}J(t)
\geq c_1(J(t))^2-c_2,~{\rm with}~c_1:=\frac{2^{4\alpha}C'_{n,\alpha}}{\omega^2_{n-1}\Gamma^2(2\alpha)},~c_2:=C''_{n,\alpha}\|\theta_0\|^2_{L^\infty}.
\end{split}
\end{eqnarray}
To continue, we denote by $I(t)$ the solution of the following ordinary differential equation
\begin{equation*}
\left\{\ba
&\frac{d}{dt}I(t)=c_1(I(t))^2-c_2,\\
&I(0)=J(0), \ea\ \right.
\end{equation*}
which, was explicitly computed in \cite{[Jarrin-Vegara-Hermosilla]} as
\begin{align*}
I(t)=\sqrt{\frac{c_2}{c_1}}\frac{\Big(J(0)+\sqrt{\frac{c_2}{c_1}}\Big)+\Big(J(0)-\sqrt{\frac{c_2}{c_1}}\Big)e^{2\sqrt{c_1c_2}t}}{\Big(J(0)+\sqrt{\frac{c_2}{c_1}}\Big)-\Big(J(0)-\sqrt{\frac{c_2}{c_1}}\Big)e^{2\sqrt{c_1c_2}t}}.
\end{align*}
Note that $I(t)$ blows up at the finite time
\begin{eqnarray*}
T_{0}:=\frac{1}{2\sqrt{c_1c_2}}\log\frac{\sqrt{c_1}J(0)+\sqrt{c_2}}{\sqrt{c_1}J(0)-\sqrt{c_2}}>0,
\end{eqnarray*}
provided
\begin{eqnarray*}
J(0)=
\int_{\R^n}\frac{\theta_0(0)-\theta_0(x)}{|x|^n}e^{-|x|}dx
>\sqrt{\frac{c_2}{c_1}}=\frac{\omega_{n-1}\Gamma(2\alpha)}{2^{2\alpha}}\sqrt{\frac{C''_{n,\alpha}}{C'_{n,\alpha}}}\|\theta_0\|_{L^\infty}.
\end{eqnarray*}
Finally, by choosing $A(n,\alpha)=\frac{\omega_{n-1}\Gamma(2\alpha)}{2^{2\alpha}}\sqrt{\frac{C''_{n,\alpha}}{C'_{n,\alpha}}}$, \eqref{Class-initial-data}, \eqref{ordinary-differential-inequality} and the comparison principle of ordinary differential equation, it follows that $J(t)\geq I(t)$, and the functional $J(t)$ also blows up at the finite time $T_{0}>0$, leading to the desired contradiction.
The proof of Theorem \ref{the-1} is then completed.
\hfill\hfill$\square$\vskip12pt
\begin{remark}
We note that the blow-up result can be extended to equations with a fractional dissipation.
By using Proposition \ref{nonlinear-weighted-inequality-exponential} and a corresponding weighted inequality for the dissipative term, and following the proof of Theorem \ref{the-1},
we can show the finite-time blowup of smooth solutions to \eqref{dissipative} with $\gamma\in(0,\alpha)$ for radial initial data fulfilling \eqref{Class-initial-data}.
One can refer to \cite{[Zhang]} for more details.
\end{remark}


{\bf Acknowledgements.}
The author is greatly indebted to the anonymous referees for their valuable comments and suggestions.
W. Zhang was supported by the National Natural Science Foundation of China (Grant No. 12501301)
and a Science and Technology Research Project of Department of Education of Jiangxi Province, China (Grant No. GJJ2200363).



\end{document}